\documentclass[a4paper,9pt]{article}
\usepackage{amsfonts}
\usepackage{amssymb}
\usepackage{latexsym}
\usepackage{amsmath}
\usepackage{amscd}
\usepackage{amsthm}
\usepackage{graphicx}
\usepackage{epsfig}
\usepackage{hyperref}
\usepackage{url}

\usepackage{graphicx}
\usepackage{pstricks}
\usepackage{pst-all}

\def\E{\mathbb{E}}
\def\E{\mathbb{E}}

\def\N{\mathbb{N}}
\def\P{\mathbb{P}}

\def\R{\mathbb{R}}

\def\Z{\mathbb{Z}}
\def\x{\mathbf{x}}

\newcommand{\ind}{\mathbf{1}}

\newtheorem{theorem}{Theorem}[section]

\newtheorem{proposition}[theorem]{Proposition}

\newtheorem{remark}[theorem]{Remark}

\usepackage{authblk}

\makeindex

\title{Collision times of random walks and applications to the Brownian web}
\author[1]{David Coupier}
\author[2]{Kumarjit Saha}
\author[3]{Anish Sarkar}
\author[4]{Viet Chi Tran}

\affil[1]{Univ. de Valenciennes, CNRS}
\affil[2]{Ashoka University}
\affil[3]{Indian Statistical Institute, Delhi}
\affil[4]{Univ. Lille, CNRS}

\begin{document}

\maketitle

\begin{abstract}
Convergence of directed forests, spanning on random subsets of lattices
or on point processes, towards the Brownian web has made the subject of
an abundant literature, a large part of which relies on a criterion proposed
by Fontes, Isopi, Newman and Ravishankar (2004). One of their convergence
condition, called (B2), states that the probability of the event that there exists
three distinct paths for a time interval of length  $t(>0)$, all starting within a
segment of length $\varepsilon$, is of small order of $\varepsilon$.
This condition is often verified by applying an FKG type correlation inequality
together with a coalescing time tail estimate for two paths.
For many models where paths have complex interactions, it is hard to
establish  FKG type inequalities. In this article, we show that
for a non-crossing path model, with certain assumptions,
a suitable upper bound on expected first collision time among three paths
can be obtained directly using Lyapunov functions. This, in turn, provides an
alternate verification of Condition (B2). We further show that in case of
independent simple symmetric one dimensional random walks or in
case of independent Brownian motions, the expected value can be
computed explicitly. We apply this alternate method of verification of (B2) to several models in
the basin of attraction of the Brownian web studied earlier in the
literature (\cite{S67}, \cite{H71}, \cite{FLT04}).
% a good upper bound of the expected collision time between three
%paths of the directed forest, i.e. the expected time of first coalescence
%between two of the three paths, can replace (B2). In the case of random
%walks of Brownian motions, the expected collision time can be computed
%explicitly. For other cases, provided the model exhibits some homogeneity
%and good Markovian properties, Lyapunov functions can be used to
%obtain the controls on the expected collision time that interest us.
\end{abstract}
%\pagebreak
	
%\chaptercontents
	
\section{Introduction}
	
The Brownian web can be described as a collection of coalescing one-dimensional
Brownian motions starting from everywhere on a space-time plane $\R^2$.
%The difficulty is to work with uncountably many Brownian motions.
The origin of this random object goes back to the work of Arratia (see \cite{A79}),
where he studied the diffusive scaling limit of coalescing simple symmetric
random walk paths starting from every point of $2\Z$  at time $0$ and
showed that it converges to a collection of coalescing Brownian motions
starting from every point of $\R$ at time $0$. The construction of a system
of coalescing Brownian motions starting from every point in space and time
was first given by T\'{o}th and Werner in \cite{TW98} where they used it to
construct what they called ``the true self-repelling motions".  In
\cite{FINR04}, Fontes et al. finally introduced a suitable topology  so
that the Brownian web can be viewed  as a random variable taking
values in a Polish space. In particular, they developed a criterion for
convergence of  rescaled forests to the Brownian web, that we will
recall later (Theorem \ref{thm:BwebConvergenceNoncrossing1}).
	
The convergence conditions include two criteria about the number
of distinct paths for a time interval of length  $t(>0)$. In particular, their
condition (B2) states that the probability of the event that there exists
three distinct paths for a time interval of length  $t(>0)$, all starting within a
segment of length $\varepsilon$,  is of small order of $\varepsilon$.
For non-crossing path models, this condition is often verified by applying
FKG type correlation inequality together with a bound on the distribution
of the coalescence time between {\em two paths}. However, FKG is a strong
property that is difficult to verify (even to expect) for many models,
especially for models where paths have complex interactions.

The purpose of this article is to propose another method of verification
of (B2) which is applicable for a non-crossing path models with certain assumptions. We show that under 
these assumptions, it is enough to have suitable estimate
 on the tail probability of the {\em first collision time among three
paths} (see Proposition \ref{prop:AltFKG2}). Using homogeneity and
Markovian  properties of the model and using a Lyapunov function type
technique, we prove a bound on the expected first collision time among
three paths (Theorem \ref{thm:AltFKGMarkov}). This allows us to use simple Markov
inequality to obtain the tail estimate and verify condition (B2) without invoking
FKG type inequalities.

Further, we show that part of this method, can be extended to
explicitly calculate the expected  first collision time among three paths in
case of independent simple symmetric random walks  and in case of
 independent Brownian motions (Theorem \ref{thm:CollTimeExp}).
Although, these results seem to be known, the only reference that we
found was \cite{SturmSwart15} [Lemma 12]. In Section \ref{sec:ApplAltB2}
we apply our method to verify (B2) for two discrete drainage network models
studied earlier in the literature: Scheidegger model \cite{S67} and
Howard model \cite{H71}. Our domination condition (Condition (iii) 
of Proposition \ref{prop:AltFKG1}) is difficult to verify for forests in continuum spaces;
however,  it might be still possible to apply the core ideas of our method. 
As an example, the Poisson tree of Ferrari et al. \cite{FLT04,FFW05} is treated in Section
\ref{sec:continuous}.
	
It should be mentioned here that there are other approaches for studying
convergence to the Brownian web which do not require Conditon (B2) to check.
Schertzer et. al. provided such an approach applicable for non-crossing path
model (Theorem  6.6 of \cite{SSS17}). Slightly stronger version of this approach
appear in Theorem 26 of \cite{CSST18}. Both these approaches applicable
for non-crossing path models and do not require Markovian properties unless
a suitable bound on the distribution of the coalescing time between two
paths is available. Long before this,  Newman et. al.  developed an alternate
condition replacing (B2) (Theorem 1.4 of \cite{NRS05}),
which is applicable for crossing path models also.

\section{Convergence to the Brownian web}\label{sec:mainresult}
	
\subsection{The Brownian web}
	
We describe the Brownian web and the topology introduced in
Fontes et al. \cite{FINR04}, that we are going to work with.
This subsection ends by recalling the convergence criterion
proposed by these authors.
	
Let $\R^{2}_c$ denote the completion of the space time plane $\R^2$ with
respect to the metric
\begin{equation*}
\rho((x_1,t_1),(x_2,t_2)) := |\tanh(t_1)-\tanh(t_2)|\vee \Bigl| \frac{\tanh(x_1)}{1+|t_1|}
-\frac{\tanh(x_2)}{1+|t_2|} \Bigr|.
\end{equation*}
% I think it is good to include this line,
As a topological space $\R^{2}_c$ can be identified with the
continuous image of $[-\infty,\infty]^2$ under a map that identifies the line
$[-\infty,\infty]\times\{\infty\}$ with the point $(\ast,\infty)$, and the line
$[-\infty,\infty]\times\{-\infty\}$ with the point $(\ast,-\infty)$.
A path $\pi$ in $\R^{2}_c$ with starting time $\sigma_{\pi}\in [-\infty,\infty]$
is a mapping $\pi :[\sigma_{\pi},\infty]\rightarrow [-\infty,\infty] \cup \{ \ast \}$ such that
$\pi(\infty)= \ast$ and, when $\sigma_\pi = -\infty$, $\pi(-\infty)= \ast$.
Also $t \mapsto (\pi(t),t)$ is a continuous
map from $[\sigma_{\pi},\infty]$ to $(\R^{2}_c,\rho)$.
We then define $\Pi$ to be the space of all paths in $\R^{2}_c$
with all possible starting times in $[-\infty,\infty]$.
	
The following metric, for $\pi_1,\pi_2\in \Pi$
\begin{align*}
d_{\Pi} (\pi_1,\pi_2) := \max \Bigl \{ & |\tanh(\sigma_{\pi_1})-\tanh(\sigma_{\pi_2})|,\\
& \sup_{t \geq \sigma_{\pi_1} \wedge \sigma_{\pi_2}} \Bigl|\frac{\tanh(\pi_1(t\vee\sigma_{\pi_1}))}{1+|t|}-\frac{
\tanh(\pi_2(t\vee\sigma_{\pi_2}))}{1+|t|}\Bigr|\Bigr\}
\end{align*}
makes $\Pi$ a complete, separable metric space. The metric $d_{\Pi}$ is slightly different
from the original choice in \cite{FINR04} which is somewhat less natural as explained
in \cite{SS08}.
Convergence in this metric can be described as locally uniform convergence of paths
as well as convergence of starting times. In what follows, for $\x = (\x(1), \x(2)) \in \R^2$ the
notation $\pi^{\x}$ is used to denote a path starting from $\x$, i.e., $\pi^{\x} : [\x(2),\infty) \mapsto \R$.

Let ${\mathcal H}$ be the space of compact subsets of $(\Pi,d_{\Pi})$ equipped with
the Hausdorff metric $d_{{\mathcal H}}$ given by,
\begin{equation*}
d_{{\mathcal H}}(K_1,K_2) := \sup_{\pi_1 \in K_1} \inf_{\pi_2 \in
	K_2}d_{ \Pi} (\pi_1,\pi_2)\vee
\sup_{\pi_2 \in K_2} \inf_{\pi_1 \in K_1} d_{\Pi} (\pi_1,\pi_2).
\end{equation*}
Since $(\Pi,d_{\Pi})$ is Polish, $({\mathcal H},d_{{\mathcal H}})$ is also Polish. Let
$B_{{\mathcal H}}$ be the Borel  $\sigma-$algebra on the metric space $({\mathcal H},d_{{\mathcal H}})$.
In the following, for any $\mathbf{x}^1, \dotsc, \mathbf{x}^k \in \R^2$,
the notation $(W^{\mathbf{x}^1}, \dotsc, W^{\mathbf{x}^k})$ represents coalescing Brownian motions
starting from $\mathbf{x}^1,\dotsc, \mathbf{x}^k$ respectively.
The Brownian web ${\mathcal W}$ is characterized as (Theorem 2.1 of \cite{FINR04}):
% We now recall from \cite{FINR04} Theorem 2.1, the following characterization
% of the Brownian web ${\mathcal W}$:
\begin{theorem}
\label{theorem:Bwebcharacterisation}
There exists an $({\mathcal H}, {\mathcal B}_{{\mathcal H}})$ valued random variable
${\mathcal W}$ such that whose distribution is uniquely determined by
the following properties:
\begin{itemize}
\item[$(a)$] for each deterministic point $\mathbf{x}\in \R^2$
there is a unique path $\pi^{\mathbf{x}}\in {\mathcal W}$  almost surely;

\item[$(b)$] for a finite set of deterministic points $\mathbf{x}^1,\dotsc, \mathbf{x}^k \in \R^2$,
the collection $(\pi^{\mathbf{x}^1},\dotsc,\pi^{\mathbf{x}^k})$ is
distributed as $(W^{\mathbf{x}^1}, \dotsc, W^{\mathbf{x}^k})$;
		
\item[$(c)$] for any countable deterministic dense set ${\mathcal D}\in \R^2$,
${\mathcal W}$ is the closure of $\{\pi^{\mathbf{x}}: \mathbf{x}\in {\mathcal D} \}$ in $(\Pi, d_{\Pi})$  almost surely.
\end{itemize}
\end{theorem}
Note that the above theorem shows how to work with uncountably many
Brownian motions as the collection is
almost surely determined by only countably many coalescing Brownian motions.

\subsection{Convergence criterion of Fontes et al.}

Since the random variable ${\mathcal W}$ takes values in a Polish space,
we can talk about weak convergence. Fontes et al. \cite{FINR04}
provided the criteria to study weak convergence to
the Brownian web ${\cal W}$.  We present their convergence criteria
for models with non-crossing paths only. A subset $\Gamma$ of $\Pi$ is
said to be consists of non-crossing paths if there exist no $\pi_1,
\pi_2 \in \Gamma$ such that
\begin{equation}
\label{def:eqnPathNonCross}
(\pi_1(t) - \pi_2(t))(\pi_1(s) - \pi_2(s)) < 0 \text{ for some }t,s \geq \max\{
\sigma_{\pi_1}, \sigma_{\pi_2}\}.
\end{equation}
For a subset $\Gamma \subseteq \Pi$ of paths  and for $t\in \R$ let $
\Gamma^{t} := \{\pi\in\Gamma  : \sigma_{\pi} \leq t\}$ denote the set
of paths which start `before' time $t$.
For $t>0$ and $t_0,a,b\in \R$ with $a<b$, we define a counting
random variable as
follows
\begin{align}
\label{def:EtaPtSet}
\eta_\Gamma(t_0,t;a,b) & := \#\{\pi(t_0+t):\pi \in \Gamma^{t_0}\text{ and }\pi(t_0)\in [a,b]\}.
\end{align}
The following theorem of \cite{FINR04} provides convergence criteria
to the Brownian web for non-crossing path models.
\begin{theorem}
\label{thm:BwebConvergenceNoncrossing1}
Let ${\mathcal X}_n, \; n \in \N$ be $({\mathcal H},B_{{\mathcal H}})$
valued random variables with non-crossing paths. Assume that the following
conditions hold:
\begin{itemize}
\item[$(I_1)$] Let ${\mathcal D}$ be a deterministic countable dense set of $\R^2$.
For each $\mathbf{x} \in {\mathcal D}$, there exists $\pi_n^{\mathbf{x}} \in {\mathcal X}_n$ such that
for any finite set of points $\mathbf{x}^1, \dotsc, \mathbf{x}^k \in {\mathcal D}$,
as $n \to \infty$, we have
$(\pi^{\mathbf{x}^1}_n, \dotsc, \pi^{\mathbf{x}^k}_n) $
converges in distribution to $(W^{\mathbf{x}^1}, \dotsc, W^{\mathbf{x}^k} )$.

\item[$(B_1)$] For all $t>0$ and $a,t_0 \in \R$, $$\limsup_{n\rightarrow \infty}\sup_{(a,t_0)\in
\R^2}\P(\eta_{{\mathcal X}_n}(t_0,t;a,a+\epsilon)\geq 2)\rightarrow 0\text{ as }\epsilon\downarrow 0.$$

\item[$(B_2)$] For all $t>0$ and $a,t_0 \in \R$, $$\frac{1}{\epsilon}\limsup_{n\rightarrow
\infty}\sup_{(a,t_0)\in \R^2}\P(\eta_{{\mathcal X}_n}(t_0,t;a,a+\epsilon)\geq 3)
\rightarrow 0\text{ as }\epsilon\downarrow 0.$$

\end{itemize}
Then ${\mathcal X}_n$ converges in distribution to the standard Brownian web ${\mathcal W}$ as $n \to \infty$.
\end{theorem}

Fontes et al. \cite{FINR04} used this convergence criteria to prove that the collection of diffusively
scaled coalescing one-dimensional simple symmetric random walk paths converges to the
Brownian web. This model is described in more detail in Subsection \ref{subsec:Scheidegger}.
Since then, the convergence
of coalescing path families (crossing as well as non-crossing) to the Brownian web has
made the subject of an abundant literature
(\cite{BGS15,CFD09,CV14,CSST18,FFW05,FINR04,NRS05,RSS16,SS13,VZ17}).
We refer to \cite{SSS17} for a review.

\subsection{Main result: an alternative verification of (B2)}
As mentioned in the introduction, for non-crossing path models in order
to verify (B2) it is enough to have sufficient tail decay for the first collision
time among three paths and often it is verified using an FKG type
correlation inequality. However,  FKG property is difficult to verify for
many models where paths have complex interactions, e.g., Howard's
model \cite{H71}, discrete directed spanning forest \cite{RSS16},
directed spanning forest \cite{BB07,CT13} etc. Among these, for
Howard's model, only a partial version of the FKG property is proved
 till now \cite{CFD09}.
The main contribution of this paper is that, we propose an alternate
method for verification of (B2)
based on a bound for expected first collision time and show that
this method is applicable
for many non-crossing path models with certain Markovian properties
and homogeneity assumptions (see Proposition \ref{prop:AltFKG2} for detail).

To fix ideas let us introduce some notations.
Suppose $ \mathcal{V} $ is a locally finite random subset of
$ \R^2$ such that $\mathcal{V}\cap \{(x,s): x\in \R, s \geq t\}$ is
nonempty for all $t \in \R$. Let $ h : \R^2 \to \mathcal{V} $ is a random
map  such that for each $(x,t) \in \R^2$, almost surely we have $h(x,t)(2) > t$
where $ h(x,t)(i) $ denotes the $i$-th co-ordinate of $ h(x,t)$ for $i=1,2$.
The vertex $h(x,t)$ is the ancestor
of $(x,t)$. We can define recursively the next ancestors, for $k \geq 1$,
by $h^k(x,t) := h(h^{k-1}(x,t))$, with $h^0(x,t) = (x,t)$.
For each $(x,t) \in \R^2$, the path $\pi^{(x,t)}: [t, \infty) \to \mathbb R$
starting from $(x,t)$ is obtained by joining the successive vertices
$ h^{k-1}(x,t),h^{k}(x,t), \; k \geq 1$ linearly. It is useful to observe
that both $h(x,t)$ and $\pi^{(x,t)}$ are defined for any $(x,t) \in \R^2$.

The collection of all the paths starting from the vertices in $\mathcal{V}$
is given as ${\mathcal X} := \{\pi^{(x,t)} : (x,t) \in \mathcal{V}\}$.
For a path $\pi \in \Pi$, for $n \geq 1$ and for
normalization constants $\gamma,\sigma > 0$, the $n$-th order diffusively scaled path
is given by $\{ \pi_n (t) = \pi_n(\gamma,\sigma)(t)=
\frac{1}{\sqrt{n}\sigma} \pi\big(n\gamma t\big)\}$. Corresponding collection of scaled paths is denoted by 
${\cal X}_n := \{\pi_n^{(x,t)}: (x,t) \in \mathcal{V}\}$.
Let $\overline{{\cal X}}_n$ denote the closure of ${\cal X}_n$ in $(\Pi, d_{\Pi})$.
We assume that for each $ n \geq 1 $, $\overline{{\cal X}}_n$  is a $({\mathcal H},
B_{{\mathcal H}})$ valued random variable.
Set $-\infty <x < y < z < \infty$ and consider the three paths $\pi^{(x,0)},
\pi^{(y,0)} , \pi^{(z,0)} $ starting from the points $(x,0), (y,0)$ and $(z,0)$
respectively.  The \textit{first collision time} of these three paths is defined as
\begin{equation}
\label{def:FirstCollTime}
T_{(x,y,z)} := \inf\{t \geq 0 : \text{ either }\pi^{(y,0)}(t) = \pi^{(x,0)}(t)
\text{ or }\pi^{(z,0)}(t) = \pi^{(y,0)}(t)\}.
\end{equation}
The next  proposition is motivated from the argument of Fontes et al.
\cite{FFW05} and uses Markov inequality in place of FKG inequality. 

\begin{proposition}
\label{prop:AltFKG1} 
Suppose $(\mathcal{V}, h)$ satisfies the following conditions:
\begin{itemize}
\item[(i)]  [Non-crossing] For any $(x_1,t_1)$ and $(x_2,t_2)$ in $ \R^2$, 
the paths $\pi^{(x_1,t_1)}, \pi^{(x_2,t_2)}$ are non-crossing.

\item[(ii)] [Homogeneity] For every $ a, t_0 \in \R $ and $ t, \epsilon > 0 $,
the distribution of the counting random variable $\eta_{{\mathcal X}_n}(t_0,
t_0+t;a, a+\epsilon)$ does not depend on $a, t_0$.

\item[(iii)] For each $ t  > 0 $,   we have 
\begin{align*} 
 & \frac{1}{\epsilon} \limsup_{ n \to \infty } 
\P \Bigl( \{ \pi(nt) : \pi \in {\cal X}, \sigma_{\pi} \leq 0, \pi(0) \in [0, \epsilon \sqrt{n}]\} \\
& \qquad \not \subseteq  \{ \pi^{(j,0)}(nt) : j \in [0, \lfloor \epsilon \sqrt{n}\rfloor + 1]\cap \Z\}  \Bigr) \to 0
\end{align*}
as $ \epsilon \downarrow 0 $. 

\item[(iv)] For  $x,y, z \in \R$ with $x < y < z$ and constants $C_1, C_2 > 0$
(independent of $x,y $ and $ z$), we have
\begin{equation}
\label{cond:exp-collision-time}
\E\big(T_{(x,y,z)}\big)\leq C_1 + C_2 (y-x)(z-y).
\end{equation}
\end{itemize}
Then $\{\overline{\mathcal X}_n: n \geq 1\}$ satisfies $(B_2)$.
\end{proposition}

% Condition (i) is a \textit{stronger} condition than assuming that ${\cal X}$ consists 
% of non-crossing paths only. Condition (iii) states that, almost surely for all large $n$, 
% the point set $\{ \pi(n) : \pi \in {\cal X}, \sigma_{\pi} \leq 0, \pi(0) \in [0, \epsilon \sqrt{n}]\}$ 
% is dominated by the set obtained from the collection of paths starting from the lattice 
% points $[0, \lfloor \epsilon \sqrt{n}\rfloor + 1]\cap \Z$ observed at time $n$. Observe 
% that first collision time for cadlag paths can be defined analogously and 
% Condition (i) ensures that the first collision time of $\pi^{(x,0)}, \pi^{(y,0)}, \pi^{(z,0)}$ 
% and first collision time of their respective cadlag versions $\tilde{\pi}^{(x,0)}, 
% \tilde{\pi}^{(y,0)}, \tilde{\pi}^{(z,0)}$ are same. 

\begin{proof} We prove the result for $\sigma = \gamma = 1$, the general case being similar.
From Condition (ii) [homogeneity], we observe that, for any $ a, t_0 \in \R $ and $ t, \epsilon > 0 $,
\begin{align*}
\P(\eta_{{\mathcal X}_n}(t_0,t_0+t;a, a+\epsilon)\geq 3)   = \P(\eta_{{\mathcal X}_n}(0,t;0,\epsilon)\geq 3).
\end{align*}
Hence, it suffices to show that for any $t>0$,
\begin{align}
\label{etape2}
\frac{1}{\epsilon}\limsup_{n\rightarrow \infty}\P(\eta_{{\mathcal X}_n}(0,t;0,\epsilon)\geq 3) 
= \frac{1}{\epsilon}\limsup_{n\rightarrow \infty}\P(\eta_{{\mathcal X}}(0,nt;0, \epsilon\sqrt{n})\geq 3)
\rightarrow 0, 
\end{align}
$\text{ as } \epsilon \downarrow 0.$
%Let us consider the unscaled paths (i.e., $ {\mathcal X}$). 
Equation \eqref{etape2} means, for $\mathcal{X}$, that there exist three
distinct paths starting between $0$ and $\epsilon \sqrt{n}$ that
do not coalesce before $nt$. We prove \eqref{etape2} for $t=1$.
The argument for general $t > 0$ is similar.

Let us define the set of paths starting from vertices $ \{ (j,0) : 0 \leq j \leq  \lfloor 
\epsilon \sqrt{n}\rfloor + 1 \} $ by $ {\cal P} $. Let  $B_{n} (\epsilon) $ 
denote the event $ \bigl\{ \{ \pi(n) : \pi \in {\cal X}_1, \sigma_{\pi} \leq 0, \pi(0) \in [0, \epsilon \sqrt{n}] \} 
 \not \subseteq \{ \pi (n) : \pi \in {\cal P} \} \bigr\} $. On the complement of the event $ B_n (\epsilon)$, there must 
be at least $3$ distinct paths in  $ {\cal P} $ till time $n$,  the event which 
we denote by $ A_{{\cal P} } (3; n) $. Therefore, we have, 
\begin{align*}
\{ \eta_{{\mathcal X}_1}(0,n;0,\epsilon \sqrt{n})\geq 3 \} \subseteq 
 A_{{\cal P} } (3; n)   \cup B_n (\epsilon).
\end{align*}

Next, we estimate  $\P (  A_{{\cal P} } (3; n)  )$. 
We follow a modification of the argument of Fontes et al. \cite{FFW05}  
and use the straightforward Markov inequality. For any $ i \in \Z $ and 
$ j \geq 2 $, define the event where no coalescence occur before time 
$n$ between the paths starting from the vertices $ (i,0), (i+1,0) $ and 
$ (i+j, 0) $, i.e.,
\begin{equation}
\label{def:Eijn}
E (i, j ; n )  := \{ \pi^{(i,0)}(n) < \pi^{(i+1,0)}(n) < \pi^{(i+j, 0)}(n) \}.
\end{equation}
Note that the paths $\pi^{(i,0)}, \pi^{(i+1,0)}$ and $\pi^{(i+j,0)}$ 
are in $ {\cal P} $ but need not be in ${\cal X}$.
% 
% 
% % From Condition (iii) it follows that for all $n$, on the event $\{\eta_{{\mathcal X}_1}
% % (0,nt;0, \epsilon\sqrt{n})\geq 3\}$, there must exist three paths among the set 
% % of paths $\{\pi^{(j,0)} : j \in [0,\lfloor \epsilon \sqrt{n} \rfloor + 1 ]\}$ such that, 
% % no two of them intersect by time $n$. 
% 
From non-crossing nature of paths (Condition (i)) it follows, if $ A_{{\cal P} } (3; n) $ occurs
then at least one of $ E (i, \lfloor \epsilon\sqrt{n}\rfloor + 1 - i ; n ) $ must occur for some 
$ i = 0, 1,  \dotsc, \lfloor \epsilon\sqrt{n}\rfloor - 1 $. 
%  that the event that 
% there exist three paths in $  {\cal P} $ such that 
% %$\{\pi^{(j,0)} : j \in [0,\lfloor \epsilon \sqrt{n} \rfloor + 1 ]\}$ such that, 
% no two of them intersect by time $n$ can be expressed as 
% $\cup_{ i = 0 }^{\lfloor \epsilon\sqrt{n}\rfloor - 1} E (i, \lfloor \epsilon\sqrt{n}\rfloor + 1 - i ; n ) $.
Hence we have:
\begin{equation}
\label{def:3PathContainedEijn}
 A_{{\cal P} } (3; n) \subseteq \cup_{ i = 0 }^{\lfloor \epsilon\sqrt{n}\rfloor - 1}
E (i, \lfloor \epsilon\sqrt{n}\rfloor + 1 - i ; n ) .
\end{equation}
Using Condition (iv) and the Markov inequality,
\begin{align*}
& \P ( E (i, \lfloor \epsilon\sqrt{n}\rfloor + 1 - i ; n ) )
 =  \P( T_{(i,i+1,\lfloor \epsilon\sqrt{n}\rfloor + 1)} > n) \\
& \leq  \frac{1}{n} \E\bigl(T_{(i,i+1,\lfloor \epsilon\sqrt{n}\rfloor + 1)} \bigr)
 \leq  \frac{1}{n} \bigl[ C_1 + C_2 \big(\lfloor \epsilon\sqrt{n}\rfloor -i\big) \bigr] \\
& \leq  \frac{1}{n}  \bigl[ C_1 + C_2  \epsilon\sqrt{n} \bigr].
\end{align*}
Thus, from (\ref{def:3PathContainedEijn}),
\begin{align*}
& \limsup_{n\to\infty} \P(\eta_{{\mathcal X}_n}(0,1;0,\epsilon)\geq 3)\\
& \leq    \limsup_{n\to \infty} \Bigl[ \frac{1}{n} \sum_{i = 0}^{\lfloor \epsilon \sqrt{n} \rfloor - 1}
\bigl[ C_1 +  C_2  \epsilon\sqrt{n}  \bigr] +  \P (B_n (\epsilon))   \Bigr]\\
& \leq    \limsup_{n\to \infty} \frac{1 }{n}  \bigl[ C_1  \epsilon \sqrt{n}
+ C_2 \big(\epsilon \sqrt{n} \big)^2 \bigr]   + \limsup_{n\to \infty}\P (B_n (\epsilon))  \\
& \leq  C_2 \epsilon^2 + \limsup_{n\to \infty} \P (B_n (\epsilon)) .
\end{align*}
Hence, using condition (iii),  $\tfrac{1}{\epsilon} \limsup_{n\to\infty} \P(\eta_{{\mathcal X}_n}(0,1;0,\epsilon)\geq 3) \to 0$
as $\epsilon \downarrow 0$.  %This completes the proof.
\end{proof}

\section{First collision time of three paths}\label{sec:CoaTime}

Conditions (i), (ii) and (iii) of Propositions \ref{prop:AltFKG1}
are usually easy to read from the model. We focus in this section on  Condition (iv).
We first show that when the paths of ${\cal X}_n$ are non-crossing and exhibit certain 
Markov properties, then \eqref{cond:exp-collision-time}
can be obtained by Lyapunov techniques. This is done in Section \ref{sec:Lyapunov}.
In case of coalescing simple symmetric  one-dimensional random walks or in case of
Brownian motions, we calculate the expectation of the first collision time explicitly.
This is presented in Sections \ref{sec:randomwalks} and \ref{sec:BM}.
%When the exact computation is not possible, \eqref{cp,cond:exp-collision-time} can be obtained by Lyapunov considerations.

\subsection{Control of the collision time expectation}\label{sec:Lyapunov}
We start with a preliminary result on entrance time of a Markov chain with countable state space.
%Later we will state the analogs version for a continuous state space markov chain.

\subsubsection{An entrance time result for a Markov chain}\label{sec:entrance}
%Some of the above results may be generalized in the Markov chain set-up.
%Although exact values of the expectation can not always be computed, bounds
%on the expectation can be obtained, showing that Condition (iv) of Propositions
%\ref{prop:AltFKG1} and \ref{prop:AltFKG1Lattice} is fulfilled.
Let $\{Y_j : j \geq 0\}$ be such a Markov chain with countable state space
${\cal M}$. For any subset $M \subseteq {\cal M}$, let $\tau(M) := \inf\{n \geq 1: Y_n \in M\}$
be the first entrance time to $M$.
% We denote $\mathbf{x} = (\mathbf{x}(1),\mathbf{x}(2)) \in {\cal M}$.
%We show the following:
\begin{theorem}
\label{thm:AltFKGMarkov}
Suppose there exist a function $V : {\cal M} \to [0,\infty)$, $b \geq 0 $, $p_0 > 0$
and two disjoint subsets $M_0, M_1$ of $ {\mathcal M}$ such that
\begin{itemize}
\item[(i)] $\E[V(Y_1) - V(Y_0) | Y_0 = x] \leq -1 + b{\mathbb I} \{x \in M_1\}$
for all $x \in M_0^c$;
\item[(ii)] $\P(Y_1 \in M_0|Y_0 = x) \geq p_0$ for all $x \in M_1$
%\item[(iii)] $\sup_{x \in M_1}V(x) < \infty$.
\end{itemize}
where $ {\mathbb I} (A) $ denote the indicator of the set $A$. Then for $x \in M_0^c$ we have
\begin{equation}
\label{eqn:MainBound}
\E[\tau(M_0) | Y_0 = x] \leq  V(x) + b/p_0 .
\end{equation}
% where $ C_1 = (b + \sup_{x \in M_1}V(x))/p_0 $.
\end{theorem}

Note here that $M_1 \cup M_0$ need not be finite. In our applications, we would have
$ V(x) = 0 $ for $ x \in M_0$ and the class of states in $ M_0 $ is closed, i.e., if we enter
this class, we would not get out of it. %Hence, for $ x \in M_0$, the left hand side of condition (i)
%would be $0$, which would put an automatic restriction of lower bound $1$ for $b$.

Let us briefly explain the intuition behind this theorem.
Condition (i) is essentially a Lyapunov type condition, where we have a strictly positive drift
towards $M_1 \cup M_0$, which by standard arguments would
imply that the expected entrance time into $M_1 \cup M_0$ for any $ x \in (M_1 \cup M_0)^c$
is bounded by constant multiple of $ V(x)$.
Condition (ii) ensures that each time the chain is in $M_1 $, there is a strictly positive
probability (bounded below) to hit $M_0$ at the next step and the chain should enter the set $ M_0$ 
after a geometric number of time points. % and condition (iii) ensures that if the chain
%does not hit $M_0$, it remains in a region with bounded values of $ V(\cdot) $.
If $M_1$ is the empty set, %Conditions (ii) and (iii) hold trivially, and 
the theorem is a known Lyapunov result.
When $M_1 \neq \emptyset$, this theorem shows that it is sufficient to focus on the entrance time to $M_0 \cup M_1$.

%In the case of the collision time between three paths of a Markov forest,
%the Markov chain can be $Y=(D^{(L,M)},D^{(M,R)})$, $M_0=\{(x,y)\in
%\R^2_+ :  xy=0\}$ and $M_1=\{(x,y)\in \R^2_+ : x\leq m\mbox{ or } y\leq m \}$, for example.

\begin{proof}[Proof of Theorem \ref{thm:AltFKGMarkov}] 
% 
%
% For the ease of notation let  $ T_1 := \tau(M_0 \cup M_1)$
% and $T_0 := \tau(M_0)$. For $x, y \in {\mathcal M}$ let $\P(x,y) :=
% \P(Y_1 = y \mid Y_0 = x)$ denote the transition probability and $\P_{x}(\cdot) :=
% \P( \cdot \mid Y_0 = {x})$ and $ \E_{x} ( \cdot ) := \E ( \cdot | Y_0 = x) $ denote respectively the law
% and expectation conditionally to the initial condition.
Let us  define $ W_0 = 0 $ and for $ n \geq 1 $, 
\begin{equation}
\label{eqn:defZ_n}
W_n = V(Y_n) - V(Y_0) - \sum_{i=1}^{n} \E \bigl[ V (Y_i) - V(Y_{i-1}) \mid Y_{i-1} \bigr].
\end{equation}
Letting $ {\cal F}_n = \sigma( Y_0, Y_1, \dotsc, Y_n) $, it is easy to check that that 
$ \{ W_n : n \geq 0 \} $ is ${\cal F}_n $-adapted martingale. 

Let us set $ Y_0  = x \in M_0^c $ and  $T_0 := \tau(M_0)$. 
Clearly $ T_0$  is ${\cal F}_n $-measurable
stopping time. Therefore, we have that  $ \{ W_{n \wedge T_0} : n \geq 0 \} $ is also ${\cal F}_n $-adapted martingale. 
Let $\P_{x}(\cdot) :=  \P( \cdot \mid Y_0 = {x})$ and $ \E_{x} ( \cdot ) := \E ( \cdot | Y_0 = x) $ 
denote respectively the law  and the expectation conditionally to the initial condition.
Hence, we  have $ \E_{x} ( W_{n \wedge T_0} ) = \E_{x} (W_{0 \wedge T_0}) = 0 $ for all $ n \geq 1 $. 
Thus, using condition (i), we obtain that 
\begin{align}
\label{eqn:V_nTauBound}
 \E_{x} \bigl[ V(Y_{n \wedge T_0}) \bigr] &  = V(x) + \E_{x} \Bigl[  \sum_{i=1}^{n \wedge T_0} 
\E_{x} \bigl[ V (Y_i) - V(Y_{i-1}) \mid Y_{i-1} \bigr] \Bigr] \nonumber \\
& \leq  V(x) + \E_{x} \Bigl[  \sum_{i=1}^{n \wedge T_0} 
-1 + b{\mathbb I} \{Y_{i-1} \in M_1\}  \Bigr] \nonumber \\ 
& =  V(x) - \E_{x} (n \wedge T_0 )  + b  \E_{x} \Bigr[ \sum_{i=1}^{n \wedge T_0} 
{\mathbb I} \{Y_{i-1} \in M_1\}  \Bigr] 
\end{align}
since for any $ n \geq 1$ and $ i = 0, 1, \dotsc, n \wedge T_0 - 1$, 
$ Y_i \in M_0^c $. Hence, transferring $ \E_{x} (n \wedge T_0 )  $ to the other side and 
observing that $ V $ is a non-negative function, we have 
\begin{equation*}
 \E_{x} (n \wedge T_0 ) \leq \E_{x} (n \wedge T_0 ) + \E_{x} \bigl[ V(Y_{n \wedge T_0}) \bigr] 
\leq V(x)  + b  \E_{x} \Bigr[ \sum_{i=1}^{n \wedge T_0} 
{\mathbb I} \{Y_{i-1} \in M_1\}  \Bigr] .
\end{equation*}

We argue below that $ T_0 < \infty $ almost surely and  $ \sum_{i=1}^{n \wedge T_0} 
{\mathbb I} \{Y_{i-1} \in M_1\}  $ is stochastically dominated by a geometric random variable 
(total number of trials before a success) with success probability $ p_0 $. Hence applying MCT
on the left hand side of above equation, we obtain the required result.

% For $x \in {\cal M}$, define
% \begin{equation*}
% u_n^{({x})} := \sum_{ k = 0}^n \P_{x} (T_0 > k)  \text{ and }
% v_n^{({x})} := \sum_{ k = 0}^n \P_{x} (T_1 > k).
% \end{equation*}
% Set $ U_n := \sup \{ u_n^{({x})}  : x \in M_1 \}$. Observe that $ u_n^{(x)}
% \leq n+1 $ so that the $ U_n $ is well defined
% for each $ n \geq 1 $.

Let  $ T_1 := \tau(M_0 \cup M_1)$. From Theorem 11.3.4 in Meyn and Tweedie \cite{MT93},
assumption (i) implies that for all $x   \not\in M_0 \cup M_1$,  $ \E_x ( T_1  )  \leq V(x) $ 
so that $ T_1 < \infty $ almost surely. 
%Though this is a standard use of the Lyapunov 
% method, we include this proof for the sake of completeness. 
% Define,  $ Z_n = V(Y_n) {\mathbb I}(T_1 >n ) $ for $ n \geq 0 $. 
% For $ x \not\in M_0 \cup M_1$,  using (i),
% we have that,  on the set $ \{ T_1 > n \} $,
% \begin{equation*}
% \E_{x} [ Z_{n+1} \mid {\cal F}_n ] \leq  {\mathbb I}(T_1 > n) \E_{x} [ V(Y_{n+1} )
% \mid {\cal F}_n ] \leq Z_n - {\mathbb I}(T_1 > n)
% \end{equation*}
% Again, on $ \{T_1 \leq n\} $ and hence $ Z_n = Z_{n+1} = 0 $ a.s. Thus, we have a non-negative
% super-adapted martingale $ \{Z_n : n \geq 0 \} $ and hence
% \begin{equation*}
%  0 \leq \E_{x} [Z_{n+1} ]  \leq  \E_{x} [ Z_n ] - \P_{x} (T_1 > n) \leq \dotsb \leq \E_x [Z_0]
% - \sum_{k=0}^n \P_{x} (T_1 > k) .
% \end{equation*}
% Using that $ Z_0 = V(x)$ and letting $ n \to \infty $, we have, for $x \not\in  M_0\cup M_1$,
% \begin{equation*}
% \sum_{k=0}^n \P_{x} (T_1 > k)   \uparrow \E_x ( T_1  )  \leq V(x).
% \end{equation*}
Now, for $ x \in (M_0 \cup M_1)^c $, the chain will almost 
surely enter the set $ M_0 \cup M_1 $ in finite time and when it enters $ M_0 \cup M_1 $, 
if it hits the set $ M_0 $ we are done. If it hits the set $ M_1 $, it has a strictly positive probability
$ p_0 $ of hitting the set $M_0$ at the next step. So, in finitely many steps, the chain will either 
hit $ M_0 $ or will go out to $ (M_0\cup M_1 )^c $. Again, if the chain hits $ M_0 $, we are done. If the 
chain is in $ (M_0\cup M_1 )^c $, it is back at the starting situation. So, again it will hit the set $   M_0 \cup M_1 $ 
in finite time and so on. Since $ p_0 > 0 $, the chain can only go out of $ M_1 $ to $ (M_0\cup M_1 )^c $ 
finitely many times, before it hits the set $ M_0 $. Thus, the chain will hit the set $ M_0 $ in finite time 
almost surely. When the chain starts from $ M_1 $, the situation is as above without first having to hit 
the set $ M_0 \cup M_1 $. So, all cases, $ T_0 < \infty $ almost surely. 

For every time point $i < n \wedge T_0 $, the chain spends in $ M_1 $ it has a probability of at least 
$ p_0 $ of hitting the set $M_0$. Thus, taking hitting $ M_0 $ at the next time point as success, it is like 
conducting a geometric trial with success probability at least $ p_0 $ and hence the number of time 
points $i < n \wedge T_0 $ is stochastically dominated by a geometric random variable with success 
probability $ p_0 $. For a formal proof, we will modify the chain's  behaviour by making each state 
in $M_0$ an absorbing state. Since we are only interested in time points before the chain 
visits $M_0$, the distribution of  $ \sum_{i=1}^{n \wedge T_0}  {\mathbb I} \{Y_{i-1} \in M_1\} $ 
is unaffected by this modification. 
With a little abuse of notation, we will continue to denote the modified chain by $ Y_n $. 
Next define $ \sigma_0 = 0 $ and  for $ i \geq  1 $, $ \sigma_{i} = 
\inf \{ n > \sigma_{i-1} : Y_i \in M_1 \} $.  Here we follow the convention that the infimum of an 
empty set $ + \infty $. Observe that $ \sigma_i $ is a $ {\cal F}_n $ 
stopping time.  Let $ N = \sup \{ k : \sigma_k < \infty \} $. We also observe that 
\begin{equation}
 \sum_{i=1}^{n \wedge T_0}  {\mathbb I} \{Y_{i-1} \in M_1\}  \leq  \sum_{i=1}^{ \infty } 
{\mathbb I} \{Y_{i-1} \in M_1\}  = N . 
\end{equation}

Now, using the fact that each state in $ M_0 $ is absorbing,  we have for $ k \geq 1 $, 
\begin{align*}
\P_x ( N > k) & = \P_x ( \sigma_{k+1} < \infty ) \\
& = \P_x ( Y_{\sigma_i} \in M_1, Y_{\sigma_i + 1} \in M_0^c \text{ for } 1 \leq i \leq k, \sigma_{k+1} < \infty ) \\
& \leq \P_x ( Y_{\sigma_i} \in M_1, Y_{\sigma_i + 1} \in M_0^c \text{ for } 1 \leq i \leq k ) \\
& = \E_{x} \Bigl[ \E \Bigl( \prod_{i=1}^k \mathbb{I} (Y_{\sigma_i} \in M_1, Y_{\sigma_i + 1} \in M_0^c) 
  \mid {\cal F}_{\sigma_k} \Bigr) \Bigr]\\
& = \E_{x} \Bigl[ \prod_{i=1}^{k-1} \mathbb{I} (Y_{\sigma_i} \in M_1, Y_{\sigma_i + 1} \in M_0^c) 
 \mathbb{I} (Y_{\sigma_k} \in M_1 ) 
 \P_{Y_{\sigma_k} } \bigl(  Y_1 \in  M_0^c \bigr)  \Bigr]\\
& \leq (1-p_0) \E_{x} \Bigl[ \prod_{i=1}^{k-1} \mathbb{I} (Y_{\sigma_i} \in M_1, Y_{\sigma_i + 1} \in M_0^c) \Bigr] \\
& \leq (1-p_0)^k
\end{align*}
following the above conditioning argument step by step. This proves the stochastic domination. 

\end{proof}

\subsubsection{Application to the expectation of first collision times }
Let us construct a discrete drainage network model where the  vertex set
${\cal V}$ is a random subset of $\Z^2$ satisfying some homogeneity
and Markovian properties (as mentioned in Proposition \ref{prop:AltFKG2})
and apply Theorem \ref{thm:AltFKGMarkov} to show that the Condition (iii)
of Proposition \ref{prop:AltFKG1} holds. We assume that for any $(x,t)\in \Z^2$,
 it's ancestor $h(x,t)$ sits on the next level, i.e., $h(x,t)(2) = t+1$ where
$ h(x,t)(i) $ is the $i$-th co-ordinate of $ h(x,t)$ for $ i =1,2$.
As in Proposition \ref{prop:AltFKG1}, $\mathcal{X}$ denotes the collection
of all paths starting from ${\cal V}$ and $\mathcal{X}_n$ is the corresponding
collection of $n$-th order diffusively scaled paths.
Define $\mathcal{S} = \{(u,v)\in \Z^2 : u, v \geq 0\}$ and $\mathcal{S}_0=
\{(u,v)\in \mathcal{S} : uv=0\}$. $\mathcal{S}$ (resp. ${\cal S}_0$) plays
the role of $\mathcal{M}$ (resp. $\mathcal{M}_0$) in Theorem \ref{thm:AltFKGMarkov}.
We are now ready to state our Proposition regarding alternate verification of (B2).

\begin{proposition}
\label{prop:AltFKG2}
Suppose $(\mathcal{V}, h)$ satisfies the following conditions:
\begin{itemize}
\item[(i)]  [Markov property] For $x,y,z \in \Z$ with $x< y<z$,
the process $\{  X_j   := ({\pi}^{(y,0)}(j) - {\pi}^{(x,0)}(j),{\pi}^{(z,0)}(j)-{\pi}^{(y,0)}(j)) : j \geq 0\}$
forms a Markov chain with state space ${\mathcal S}$.

\item[(ii)] [Control of the collision time] There exist a subset
${\mathcal S}_1$ of ${\cal S} \setminus {\mathcal S}_0$ and $p_0 > 0$, such that
\begin{itemize}
\item[(a)] there exist positive constants $d_1$ and $d_2$ such that for all $(l,m) \notin {\cal S}_0$,
\begin{align}
& \E\big(X_1(1)X_1(2) - X_0(1)X_0(2) \ |\ X_0=(l,m)\big) \nonumber \\
& \qquad  \leq - d_1 +
d_2 \mathbf{1}_{\{(l,m) \in {\cal S}_1\}}\label{cond:exp-collision-time-Lyapunov}
\end{align}
where $X_n = (X_n(1),X_n(2))$.
\item[(b)] For all $(l,m) \in {\mathcal S}_1$,  $\P(X_{1} \in {\mathcal S}_0|X_0 = (l,m)) \geq p_0$.
\end{itemize}
\end{itemize}
Then for any $-\infty <x < y < z < \infty$, there exist positive constants $C_1, C_2$ (independent of $x,y,z$) such that 
$$
\E(T_{(x,y,z)}) \leq C_1 + C_2 (y-x)(z-y).
$$
In addition, if $(\mathcal{V}, h)$ satisfies Condtion (ii)  of Proposition 
\ref{prop:AltFKG1}, then $\{\overline{\mathcal X}_n: n \geq 1\}$ satisfies $(B2)$.
\end{proposition}

Note that we have considered the state space of the Markov chain as the first quadrant, which implies 
that the paths are automatically non-crossing. Since the vertex set ${\cal V}$ is a random subset of $\Z^2$, 
and for any $(x,t) \in \Z^2$ we have $h(x,t)(2) = t+1$ almost surely, hence Condition (iii) of Proposition \ref{prop:AltFKG1} automatically holds. 

\begin{proof} Consider the function $ V : {\cal S} \to [0,\infty)$ defined by
\begin{equation*}
V (u,v) = \frac{uv}{d_1}.
\end{equation*}
% We have mentioned that though Theorem \ref{thm:AltFKGMarkov} is stated 
% for countable state space Markov chain, it holds for 
% a continuous state space Markov chain as well. 
Part (a) and part (b) of 
Assumption (ii) of Proposition \ref{prop:AltFKG2} ensure that Conditions (i) and
(ii) of Theorem \ref{thm:AltFKGMarkov} hold for the markov chain $\{X_j : j \geq 0\}$. 
%Since $ {\mathcal S}_1 $ is bounded, Condition (iii) of Theorem
%\ref{thm:AltFKGMarkov} also holds. 
% Since $\{X_j : j \geq 0\}$ takes values in ${\cal S}$, first collision time 
% of the continuous paths is same as that of the cadlag versions.

For any $-\infty <x < y < z < \infty$, first collision time for the paths 
${\pi}^{(x,0)}, {\pi}^{(y,0)}$ and ${\pi}^{(z,0)}$
 is same as the entrance time to the set ${\cal S}_0$ for the Markov chain $\{X_j : j \geq 0\}$. 
Thus, the bound for  $\E(T_{(x,y,z)})$ follows from Theorem \ref{thm:AltFKGMarkov}.
Finally the proof follows from Proposition \ref{prop:AltFKG1}.
\end{proof}

\subsection{First collision time of three random walks}\label{sec:randomwalks}

%We start with a simple calculation on the first collision time of three simple symmetric random walks.
The calculations leading to equation (\ref{eqn:V_nTauBound}) yields a surprising exact computational result for
collision time among three simple symmetric random walks. With $ V(\cdot) $ being defined as the product
of the distances between the pair of the left and the middle random walk and the pair of
middle and right random walk (see below for exact definition), it turns out
that the condition (i) in Theorem \ref{thm:AltFKGMarkov} is an equality with $ M_1 $ being the empty set.
Though the exact computation is not required for the convergence to the Brownian web, we give it here
for the sake of completeness. 

Let us consider three symmetric nearest neighbour independent random walks, which we denote by
$ S_n^{(L)}, S_n^{(M)} $ and $ S_n^{(R)} $, starting from $ -2i, 0 $ and $ 2j$ respectively with $ i, j \geq 1 $.
Given three independent sequences of {\bf i.i.d}  increment random variables $ \{ I_k^{(L)} : k \geq 1 \},
\{ I_k^{(M)} : k \geq 1 \}  $ and $ \{ I_k^{(R)} : k \geq 1 \} $ with
Rademacher distributions (each of these random variables take values $ +1 $ and $ - 1$
with probability $ 1/2 $), the random walks are represented by
\[ S_n^{(L)} = -2i + \sum_{k = 1}^n I_k^{(L)},\quad   S_n^{(M)} = \sum_{k = 1}^n
I_k^{(M)} \quad \mbox{ and }\quad S_n^{(R)} = 2j + \sum_{k= 1}^n I_k^{(R)} .\]
Let us define the associated filtration $\{\mathcal{F}_n : n\geq 0\}$ 
where ${\cal F}_n = \sigma ( I_k^{(L)}, I_k^{(M)}, I_k^{(R)} :
k \leq n) $ is the $\sigma$-algebra generated by the increment random variables up to time $n$.

We define the meeting times of these random walks by
\begin{align*}
\tau_{(L,M)} & := \inf \{ n \geq 1 : S_n^{(L)} = S_n^{(M)} \}  \quad \mbox{ and }\\
\tau_{(M,R)}  &:= \inf \{ n \geq 1 : S_n^{(M)} = S_n^{(R)} \}.
\end{align*}
We know that $ \tau_{(L,M)}$ and $\tau_{(M,R)}$ are almost surely finite, but that their
expectations are infinite. Let us define the first collision time of these
random walks by
\begin{equation}
\label{def:CollLM}
\tau_C :=  \min \{ \tau_{(L,M)}, \tau_{(M,R)} \}.
\end{equation}
\begin{theorem}
\label{thm:CollTimeExp}
With the notation as above, we have
$\E (\tau_{C}) = 4ij $.
\end{theorem}

Although this theorem seems to be known, the only reference that we
 found was \cite[Lemma 12]{SturmSwart15}. We give the proof of the theorem for the sake of completeness.
Before we proceed further, let us set up the notation.
For all $n\geq 0$, we denote the distances between the paths by
\[ D_n^{(L,M)} = S_n^{(M)} - S_n^{(L)} \quad \mbox{ and }\quad D_n^{(M,R)} = S_n^{(R)} - S_n^{(M)}. \]We have $ D_0^{(L,M)} = 2i $ and $ D_0^{(M,R)} = 2j $. Also
set the maximum distance between the pairs as \[ M_n = \max \{ D_n^{(L,M)}, D_n^{(M,R)} \} .\]

\begin{remark}
It is obvious that at the first collision time, only one pair of random walks can meet.
Assuming Theorem \ref{thm:CollTimeExp}, we can also compute that the other
pair of paths is at an expected distance of $ 2(i+j)$ at the collision time.
Since at $ \tau_{C} $, at least one of $ D_{\tau_{C}}^{(L,M)}$ or
$ D_{\tau_{C}}^{(M,R)} $ is zero, we have $ M_{\tau_{C}} = D_{\tau_{C}}^{(L,M)} + D_{\tau_{C}}^{(M,R)} $. But
we may write
\begin{equation*}
D_{n}^{(L,M)} + D_{n}^{(M,R)} - 2(i+j) = \sum_{ j = 1}^n (I_j^{(R)} - I_j^{(L)} ) .
\end{equation*}
Therefore, using Wald's identity, we have
\begin{equation*}
\E ( M_{\tau_{C}} - 2(i+j) ) = \E (  D_{\tau_{C}}^{(L,M)} + D_{\tau_{C}}^{(M,R)}  - 2 (i+j)) = 0.
\end{equation*}
Thus, $ \E ( M_{\tau_{C} }) =  2(i+j). $
\end{remark}

\bigskip

The main observation is that we may rewrite $ \tau_{C} $ as
\begin{equation}
\label{def:TauDefDnMultZero}
\tau_{C} = \inf \{ n \geq 1 : D_n^{(L,M)} D_n^{(M,R)} = 0 \}.
\end{equation}

Next we have two straightforward Propositions which involve the product of the distances.

\begin{proposition}
\label{prop:ProdDistAddTimeMart}
The family $ \{ D_n^{(L,M)} D_n^{(M,R)} + n : n \geq 0 \} $ is an $ {\cal F}_n$-adapted martingale.
\end{proposition}

\begin{proposition}
\label{prop:ProdDistsAllPairsMart}
The family  $ \{ D_n^{(L,M)} D_n^{(M,R)} (D_n^{(L,M)} + D_n^{(M,R)}) : n \geq 0 \} $
is an $ {\cal F}_n$-adapted martingale.
\end{proposition}

Assuming the Propositions, we can now prove Theorem \ref{thm:CollTimeExp}.

\begin{proof}[Proof of Theorem \ref{thm:CollTimeExp}] By Proposition \ref{prop:ProdDistAddTimeMart}, we have
that, $ \{ D_{n \wedge \tau_{C}}^{(L,M)} D_{n \wedge \tau_{C}}^{(M,R)} + n \wedge \tau_{C} : n \geq 0 \} $
is a martingale and hence for any $ n \geq 1 $,
\begin{align}
\label{eqn:mainRelD_nTau}
\E ( D_{n \wedge \tau_{C}}^{(L,M)} D_{n \wedge \tau_{C}}^{(M,R)} + n \wedge \tau_{C} ) & =
\E ( D_{0 \wedge \tau_{C}}^{(L,M)} D_{0 \wedge \tau_{C}}^{(M,R)} + 0 \wedge \tau_{C} )\nonumber \\
& = \E ( D_{0}^{(L,M)} D_{0}^{(M,R)} ) = 4ij.
\end{align}
This is the version of equation (\ref{eqn:V_nTauBound}) in this case. 
Since $ \tau_{(L,M)}, \tau_{(M,R)} < +\infty $ almost surely, we have that
$ \tau_{C} < +\infty $ almost surely. By the monotone convergence theorem, we have that,
as $ n \to \infty $,
$\E ( n \wedge \tau_{C} ) \to \E(\tau_{C}) $. Since $ \tau_{C} < +\infty $ almost  surely, 
we have that, almost surely, as $n \to \infty$,
$D_{n \wedge \tau_{C}}^{(L,M)} D_{n \wedge \tau_{C}}^{(M,R)}  \to  
D_{\tau_{C}}^{(L,M)} D_{\tau_{C}}^{(M,R)} = 0$.
To complete the proof, we show that, as $ n \to \infty $,
\begin{equation}
\label{eqn:DnMinTauConv}
\E ( D_{n \wedge \tau_{C}}^{(L,M)} D_{n \wedge \tau_{C}}^{(M,R)} ) \to 0.
\end{equation}

From Proposition \ref{prop:ProdDistsAllPairsMart}, we have that $ \{  D_{n \wedge \tau_{C}}^{(L,M)}
D_{n \wedge \tau_{C}}^{(M,R)} (D_{n \wedge \tau_{C}}^{(L,M)} + D_{n \wedge \tau_{C}}^{(M,R)} ) : n \geq 0 \} $
is also a martingale and hence,  for any $ n \geq 1 $,
\begin{align}
& \E \bigl[  D_{n \wedge \tau_{C}}^{(L,M)}
D_{n \wedge \tau_{C}}^{(M,R)} (D_{n \wedge \tau_{C}}^{(L,M)} + D_{n \wedge \tau_{C}}^{(M,R)} ) \bigr]\nonumber \\
& =   \E \bigl[  D_{0 \wedge \tau_{C}}^{(L,M)}
D_{0 \wedge \tau_{C}}^{(M,R)} (D_{0 \wedge \tau_{C}}^{(L,M)} + D_{0 \wedge \tau_{C}}^{(M,R)} )  \bigr]\nonumber\\
& =  \E \bigl[  D_{0 }^{(L,M)}  D_{0}^{(M,R)} (D_{0}^{(L,M)} + D_{0 }^{(M,R)} )  \bigr] = 8ij (i+j).\label{eqn:main2RelD_nTau}
\end{align}
Therefore, observing that $ D_{n \wedge \tau_{C}}^{(L,M)}  $ and $  D_{n \wedge \tau_{C}}^{(M,R)}  $
are both non-negative, we have
\begin{align*}
& \sup_{ n \geq 1 } \E \bigl[  \bigl( D_{n \wedge \tau_{C}}^{(L,M)} D_{n \wedge \tau_{C}}^{(M,R)}  \bigr)^{3/2} \bigr] \\
& \leq  \tfrac{1}{2} \sup_{ n \geq 1 }  \E \bigl[  D_{n \wedge \tau_{C}}^{(L,M)} D_{n \wedge \tau_{C}}^{(M,R)}
(D_{n \wedge \tau_{C}}^{(L,M)}  + D_{n \wedge \tau_{C}}^{(M,R)}  ) \bigr] = 4ij(i+j)
\end{align*}
and hence $ \bigl\{ D_{n \wedge \tau_{C}}^{(L,M)} D_{n \wedge \tau_{C}}^{(M,R)}  : n \geq 1 \bigr\} $
is a sequence of uniformly integrable random variables. We conclude 
\eqref{eqn:DnMinTauConv} using Theorem 16.13. of Billingsley \cite{billingsley}.
\end{proof}

Next we prove the Propositions which are straightforward.

\begin{proof}[Proof of Proposition \ref{prop:ProdDistAddTimeMart}] This proposition
follows from straightforward calculation:
\begin{multline*}
\E \bigl[  (D_{n+1}^{(L,M)} D_{n+1}^{(M,R)} + n+1) - ( D_{n}^{(L,M)}
D_{n}^{(M,R)} +    n)  \mid {\cal F}_n \bigr]  \\
\begin{aligned}
& =   D_{n}^{(L,M)} \E \bigl[   (I_{n+1}^{(R)} - I_{n+1}^{(M)})  \mid {\cal F}_n \bigr]  +
D_{n}^{(M,R)} \E \bigl[    (I_{n+1}^{(M)} - I_{n+1}^{(L)})   \mid {\cal F}_n \bigr]  \\
&  \hspace{0.5cm} +  \E \bigl[   (I_{n+1}^{(M)} - I_{n+1}^{(L)}) ( I_{n+1}^{(R)}
- I_{n+1}^{(M)})  \mid {\cal F}_n \bigr] +1 \\
& =    \E \bigl[   I_{n+1}^{(M)} I_{n+1}^{(R)}  - ( I_{n+1}^{(M)})^2
- I_{n+1}^{(L)}  I_{n+1}^{(R)}
+  I_{n+1}^{(L)} I_{n+1}^{(M)})   \bigr] + 1 = 0.
\end{aligned}
\end{multline*}
\end{proof}

\begin{proof}[Proof of Proposition \ref{prop:ProdDistsAllPairsMart}]
Again straightforward calculations yield the result.
We have
\begin{align*}
& D_{n+1}^{(L,M)} D_{n+1}^{(M,R)} (D_{n+1}^{(L,M)} + D_{n+1}^{(M,R)} )  -
D_{n}^{(L,M)} D_{n}^{(M,R)}  (D_{n}^{(L,M)} + D_{n}^{(M,R)}) \\
& = D_{n}^{(L,M)} D_{n}^{(M,R)}  (I_{n+1}^{(R)} - I_{n+1}^{(L)} ) +
D_{n}^{(L,M)}  (D_{n}^{(L,M)} +  D_n^{(M,R)})     (I_{n+1}^{(R)} - I_{n+1}^{(M)})  \\
& \qquad + D_{n}^{(L,M)}    (I_{n+1}^{(R)} - I_{n+1}^{(M)}) (I_{n+1}^{(R)} - I_{n+1}^{(L)} )\\
& \qquad +  D_{n}^{(M,R)}  (D_{n}^{(L,M)} +  D_n^{(M,R)})   (I_{n+1}^{(M)} - I_{n+1}^{(L)}) \\
&  \qquad + D_n^{(M,R)}   (I_{n+1}^{(M)} - I_{n+1}^{(L)})  (I_{n+1}^{(R)} - I_{n+1}^{(L)} )\\
& \qquad +  (D_{n}^{(L,M)} +  D_n^{(M,R)})  (I_{n+1}^{(M)} - I_{n+1}^{(L)}) ( I_{n+1}^{(R)}
- I_{n+1}^{(M)})  \\
&  \qquad + (I_{n+1}^{(M)} - I_{n+1}^{(L)}) ( I_{n+1}^{(R)} - I_{n+1}^{(M)})
(I_{n+1}^{(R)} - I_{n+1}^{(L)}) .
\end{align*}
The conditional expectation with respect to ${\cal F}_n $ of the first term, second term and the fourth term is clearly $0$.
We observe that that the conditional expectation with respect to ${\cal F}_n $ of the  third term is $
D_{n}^{(L,M)},   $   the fifth  term is $ D_{n}^{(M,R)}   $ and the sixth term is $ -  (D_{n}^{(L,M)} +  D_n^{(M,R)})  $.
The last term is independent of $ {\cal F}_n $ and has expectation $0$. This proves the result.
\end{proof}

\subsection{First collision time of three Brownian motions}\label{sec:BM}

As one can expect, these results generalize to the Brownian motion
set up in a straightforward way. Let us start with three independent Brownian motions
$\{(B^{(L)}_t, B^{(M)}_t, B^{(R)}_t) : t\geq 0\}$ started from $ -x $, $0$ and $ y $
respectively, for $ x, y > 0 $. We also define the associated filtration $\{\mathcal{F}_t
: t\geq 0\}$ where $ {\cal F}_{t} = \sigma ( B_s^{(L)}, B_s^{(M)}, B_s^{(R)} : s \leq t ) $.

Again, let us define
\begin{equation}
\label{def:CollTimeBM}
\tau_{C}^{(B)} = \inf \{ t > 0 : (B_t^{(M)} - B_t^{(L)})(B_t^{(R)} - B_t^{(M)}) =  0 \} .
\end{equation}

Then, we have the following result
\begin{theorem}
\label{thm:CollTimeExpBM}
With the notation as above, we have
$\E (\tau_{C}^{(B)}) = xy $.
\end{theorem}

The proof follows from the propositions below which are exact replica of
Propositions \ref{prop:ProdDistAddTimeMart} and \ref{prop:ProdDistsAllPairsMart}, exactly as in the discrete case.

\begin{proposition}
We have $ \{ (B_t^{(M)} - B_t^{(L)})(B_t^{(R)} - B_t^{(M)}) + t  : t \geq 0 \} $ is a $(\mathcal{F}_t)$-adapted martingale.
\end{proposition}

\begin{proposition}
We have $ \{ (B_t^{(M)} - B_t^{(L)})(B_t^{(R)} - B_t^{(M)}) (B_t^{(R)} - B_t^{(L)})
: t \geq 0 \} $ is a $(\mathcal{F}_t)$-adapted martingale.
%adapted to the $\sigma$-algebra $ {\cal F}_{t} $ defined above.
\end{proposition}

\begin{remark}
It should be observed that there is nothing surprising in the above
Propositions as well as Propositions
\ref{prop:ProdDistAddTimeMart} and \ref{prop:ProdDistsAllPairsMart}. For example, we can write
\begin{align*}
& (B_t^{(M)} - B_t^{(L)})(B_t^{(R)} - B_t^{(M)}) + t \\
& = B_t^{(M)} B_t^{(R)}
- B_t^{(L)} B_t^{(R)} + B_t^{(L)} B_t^{(M)} -
\bigl[ (B_t^{(M)} )^2 - t \bigr]
\end{align*}
Since Brownian motions are themselves martingales and since products of independent
martingales and sums of martingales
are martingales, it is easy to observe that the first three terms are martingales. Since
 $ t $ is the quadratic variation process of the Brownian motion, the last term is also
a martingale. With these arguments, other processes can also be shown to be martingales.
\end{remark}

\section{ Applications to drainage network models }
\label{sec:ApplAltB2}

We now apply Propositions \ref{prop:AltFKG1} and \ref{prop:AltFKG2} to some
directed forests of the literature and provide an alternate verification of $(B_2)$ for these models.

\subsection{Scheidegger's model}
\label{subsec:Scheidegger}

We start with the system of coalescing simple symmetric random
walks starting from every point on $  \Z^2_{\text{even}} := \{(x,t) : x + t \text{ even}\}$.

Let $\{ b_{(x,t)}: (x,t) \in \Z^2_{\text{even}}\}$ be an i.i.d. collection such that
\begin{align*}
b_{(x,t)} =
\begin{cases}
+ 1 & \text{ with probability } 1/2 \\
- 1 & \text{ with probability } 1/2 .
\end{cases}
\end{align*}
Define $h : \Z^2 \mapsto \Z^2_{\text{even}}$ as 
\begin{align*}
h(x,t) =
\begin{cases}
(x + b_{(x,t)}, t+1) & \text{ if } (x,t) \in \Z^2_{\text{even}} \\
(x, t+1) & \text{ otherwise } .
\end{cases}
\end{align*}
For $k \geq 1$ let $h^k(x,t) := h(h^{k-1}(x,t))$ with $h^0(x,t) = (x,t)$.
For $(x,t) \in \Z^2$ let $\pi^{(x,t)}: [t, \infty) \to \mathbb R$ be the path starting at $(x,t)$
obtained by joining the successive vertices $ h^{k-1}(x,t),h^{k}(x,t), \; k \geq 1$
by straight line segments.
Let ${\mathcal X} := \{\pi^{(x,t)} : (x,t) \in \Z^2_{\text{even}}\}$,
be the collection of all paths starting from $\Z^2_{\text{even}}$.
For $n \geq 1$ and for normalization constants $\gamma,\sigma > 0$, with slight abuse of notation
let ${\cal X}_n ={\cal X}_n(\gamma,\sigma) $ denote the collection of the scaled paths.
% Let $\overline{{\cal X}}_n$ denote the closure of ${\cal X}_n$ in $(\Pi, d_{\Pi})$.
For each $ n \geq 1 $, $\overline{{\cal X}}_n$, the closure of ${\cal X}_n$ in $(\Pi, d_{\Pi})$,
is a $({\mathcal H},B_{{\mathcal H}})$
valued random variable.

Fontes et. al \cite{FINR04} proved that $\overline{{\cal X}}_n$  converges in distribution
to the Brownian web with $\gamma = \sigma = 1$ (see Theorem 6.1 of Fontes et al. \cite{FINR04}).
%The proof involves  verification of the conditions of Theorem \ref{thm:BwebConvergenceNoncrossing1}.
To verify the Condition $(B_2)$, they prove a version of FKG inequality for paths and
using it, along with the coalescing time tail estimate for a pair of random walks, to derive
suitable bounds. In this case it is not hard to obtain FKG inequalities
as the paths are independent till the time of coalescence. On the other hand, we proved in Section
\ref{sec:CoaTime} that the expected first collision time
of three simple symmetric random walks starting at $ 2i $ and $ 2j$ is exactly $ 4 ij $,
which gives the same order bound (with known constant values) as obtained by
Fontes et al. However, as a warmup, we show that Proposition
\ref{prop:AltFKG2} is applicable here.
Note that the exact value of this expectation is not required for our calculations.

In this case, we have $ \mathcal{V} = \Z^2_{\text{even}} $
and $ h(\cdot) $ is  defined as above. From the construction it is clear that the paths are non-crossing and Condition
(i) of Proposition \ref{prop:AltFKG1} is satisfied. Because of the i.i.d. nature of the collection of random
variables $ \{ b_{(x,t)} : (x,t) \in \Z^2_{\text{even}} \} $,
the homogeneity condition, i.e., Condition (ii) of Proposition  \ref{prop:AltFKG1} is also satisfied.

% For any $ (x,0) \in \Z^2_{\text{even}}  $,   $ h^{j+1} (x ,0) (1) =  h^{j} (x ,0) (1)  +
% b_{(h^{j} (x ,0) (1) , j)} $. Thus, for $ x < y < z$ with $ x,y,z \in 2 \Z$, we have
% $ ( h^{j+1} (y ,0) (1) - h^{j+1} (x ,0) (1) , h^{j+1} (z ,0) (1) - h^{j+1} (y ,0) (1))
% = ( h^{j} (y ,0) (1) - h^{j} (x ,0) (1) , h^{j} (z ,0) (1) - h^{j} (y ,0) (1)) +
% ( b_{(h^{j} (y ,0) (1) , j)} - b_{(h^{j} (x ,0) (1) , j)}, b_{(h^{j} (z ,0) (1) , j)} - b_{(h^{j} (y ,0) (1) , j)}) $.
For any $ (x,0) \in \Z^2_{\text{even}}  $,   $ h^{j+1} (x ,0) (1) =  h^{j} (x ,0) (1)  +
b_{(h^{j} (x ,0))} $. Set $ x,y,z \in 2 \Z$
such that $ x < y < z$. Then we have
$ ( h^{j+1} (y ,0) (1) - h^{j+1} (x ,0) (1) , h^{j+1} (z ,0) (1) - h^{j+1} (y ,0) (1))
= ( h^{j} (y ,0) (1) - h^{j} (x ,0) (1) , h^{j} (z ,0) (1) - h^{j} (y ,0) (1)) +
( b_{(h^{j} (y ,0))} - b_{(h^{j} (x ,0) )}, b_{(h^{j} (z ,0))} - b_{(h^{j} (y ,0))}) $.
By the independence of the family of random variables $ \{ b_{(x,t)} : (x,t) \in \Z^2_{\text{even}} \} $,
it is clear that 
\begin{align*} & \{   X_j = ({\pi}^{ (y ,0)}(j) - {\pi}^{ (x ,0)}(j), {\pi}^{ (z ,0)}(j) - {\pi}^{ (y ,0)}(j)) \\
& \qquad \qquad  = ( h^{j} (y ,0) (1) - h^{j} (x ,0) (1) , h^{j} (z ,0) (1) - h^{j} (y ,0) (1))  :
j \geq 0 \} 
\end{align*}
 is a Markov chain with state space $\{(x,t) \in \Z^2 : x,t \geq 0\}$, which 
verifies Condition (i) of Proposition  \ref{prop:AltFKG2}.

Finally, to verify Condition (ii), we take $ {\cal S}_1 = \emptyset $ so that part (b) of Condition (ii) is trivially satisfied.
For Condition (a) of Proposition \ref{prop:AltFKG2} (ii), let us denote, $b_{(x,0)} = I_x, b_{(y,0)} = I_y$
and $b_{(z,0)} = I_z$. Hence,
\begin{align*}
& \E[X_{1}(1)X_{1}(2) - X_{0}(1)X_{0}(2) | X_0 = (y-x,z-y)]  \\
& =  \E \bigl[ ( y-x + (I_y - I_x))  (z-y + (I_z - I_y))  \\
&  \qquad - (y-x)(z-y) | X_0 = (y-x,z-y) \bigr] \\
& =  \E \bigl[ ( y-x + (I_y - I_x))  (z-y + (I_z - I_y))  - (y-x)(z-y) \bigr] \\
& = (z-y) \E (I_y - I_x)) +  (y-x) \E (I_z - I_y) + \E [ (I_y - I_x)(I_z - I_y) ].
\end{align*}
Clearly from definition,  we have
$\E(I_x) =  \E(I_y) = \E(I_z) = 0$ and $ \E( I_y^2) = 1 $. Further for $x < y < z$, the random
variables $ I_x, I_y $ and $ I_z$ are independent. Hence,
\begin{equation*}
\E[X_{1}(1)X_{1}(2) - X_{0}(1)X_{0}(2) | X_0 = (y-x,z-y)] = - 1.
\end{equation*}
If $ x = y $ or $ y = z $, i.e., at least one pair of paths have already coalesced,
we have that $  X_{1}(1)X_{1}(2) = 0 $. Therefore,  ${\cal X}$ satisfies the conditions
in Proposition \ref{prop:AltFKG2} with the choice of
${\cal S}_1 = \emptyset$, $d_1 = -1$. 
%and $d_2 = 1$ (see comments after Theorem \ref{thm:AltFKGMarkov}).
This allows to conclude that  Condition (B2) of Theorem
\ref{thm:BwebConvergenceNoncrossing1} holds for $\{\overline{{\cal X}}_n : n \in \N\}$.

\subsection{Howard's model}
\label{subsec:Howard}

We first describe a 2-dimensional drainage network
proposed by Howard (see \cite{H71}).
Fix $0< p < 1$ and let $\{B_{(x,t)} :
(x,t) \in \Z^2 \}$ be
an i.i.d. collection of Bernoulli random variables with success probability $p$.
Set $ \mathcal{V} = \{ (x,t) \in \Z^2 : B_{(x,t)} = 1 \} $.
Let $\{U_{(x,t)} : (x,t) \in \Z^2 \}$ be
an i.i.d. collection of Rademacher random variables, independent of the collection
$ \{ B_{(x,t)} : (x,t) \in \Z^2 \}$: $ \P ( U_{(x,t)} = 1 ) =  \P ( U_{(x,t)} = -1)  = 1/2 $.
For a vertex  $ (x,t) \in \Z^2 $, we consider   $ k_0 =
\min \bigl\{ | k |  : k \in \Z,  B_{(x + k , t+1)} = 1 \bigr\} $. Clearly, $ k_0 $ is almost surely finite.
Now, we define,
\begin{align}
\label{eq:h(x,t)Howard}
h (x,t) :=
\begin{cases}
(x+k_0, t+1) & \text{ if }  (x - k_0 , t+1) \not\in \mathcal{V}\\
(x-k_0, t+1) & \text{ if } (x+k_0, t+1) \not \in \mathcal{V} \\
(x+U_{(x,t)}k_0, t+1) & \text { otherwise}.
\end{cases}
\end{align}
Here each open vertex $ (x,t) \in \mathcal{V}$ represents a water source and
the edge $\langle (x,t), h(x,t)\rangle$ represents the channel through
which water can flow.
Consider the random graph $G= (\mathcal{V},E)$ where the random
edge set $E$ is given by $E = \{\langle (x,t), h(x,t)\rangle : (x,t) \in \mathcal{V}\}$.
Gangopadhyay et al. \cite{GRS04} proved that this random graph $G$ is
connected almost surely. This result, shows that the  Howard's model generates
a random directed tree on $\Z^2$. The authors of \cite{GRS04} actually
used a general construction and studied this model for higher dimensions also. The
construction that we presented here is applicable  for $d=2$ only and
agrees with that of \cite{GRS04} for $d=2$.

Let ${\cal X}$ be the collection of all paths obtained by following the
edges and for $\gamma, \sigma > 0$, let ${\cal X}_n = {\cal X}_n(\gamma, \sigma)$
be the collection of scaled paths with normalization constants $\gamma,\sigma$.
Coletti et al. \cite{CFD09,CV14} proved that,
for $\gamma = 1$ and $\sigma = \sigma_0 (p)$ with $ \sigma_0 (p) $ as
given by $ \sigma_0 (p) = \bigl(\frac{(1-p)(2 - 2p + p^2)}{p^2(2 - p)^2}
\bigr)^{1/2}$, $\overline{{\cal X}}_n$  converges in distribution
to the Brownian web.
Again, to verify Condition $(B2)$ they proved a partial FKG inequality
for the Howard's model.
We show here that Howard's model satisfies Proposition \ref{prop:AltFKG2},
and thereby Condition $ (B2) $ is verified.

By the i.i.d. nature of the collections of random variables $ \{ B_{(x,t)} : (x,t) \in \Z^2\},
\{ U_{(x,t)} : (x,t) \in \Z^2\} $,
it follows that ${\cal X}$ satisfies Condition (ii) of Proposition  \ref{prop:AltFKG2}.
It is not difficult to see that ${\cal X}$ consists of non-crossing paths only and 
Condition (i) of Proposition  \ref{prop:AltFKG2} holds.
From definition it follows that for any $ (y,s) \in \Z^2 $,
$ h (y ,s)$ depends only on the random variables $ \{ B_{(x,s+1)} : x \in \Z\}$ and
$ U_{(y,s)}  $.
Thus from the i.i.d. nature of the random variables it follows that for $ x,y,z \in \Z$ with
$ x < y < z$, the process 
\begin{align*} & \{   X_j = ({\pi}^{ (y ,0)}(j) - {\pi}^{ (x ,0)}(j), {\pi}^{ (z ,0)}(j) - {\pi}^{ (y ,0)}(j)) \\
& \qquad \qquad  = ( h^{j} (y ,0) (1) - h^{j} (x ,0) (1) , h^{j} (z ,0) (1) - h^{j} (y ,0) (1))  :
j \geq 0 \} 
\end{align*}
 is a Markov chain with state space $\{(x,t) \in \Z^2: x,t \geq 0\}$, which verifies Condition (i)
of Proposition  \ref{prop:AltFKG2}.

To verify Condition (ii), we take $ \mathcal{S}_1 = \{(l,m) \in \Z^2:  1 \leq \min ( l,m ) \leq  r_0\} $
where $r_0 \in \N$ will be chosen later. We observe that, for any $ r_0 \geq 1 $,
$$\P[X_1 \in \mathcal{S}_0 \mid X_0 = (y-x,z-y)] \geq (1-p)^{2r_0} p $$
 for all $ ((y-x),(z-y)) \in {\cal S}_1 $. For this assume that $ y-x \leq r_0$ and 
 consider the event $ \{ B_{ (u,1) } = 0, $ for
$ |u-x| \leq r_0, u \neq y, B_{ (y,1) } = 1 \}$, i.e.,   the vertex above the second path is
open and all other points closed which are within a distance $r_0 $ (on both sides) from the vertex above
the first path. This ensures that the first path and the second path coalesce. Similar arguments 
hold when $ z-y \leq r_0 $. 
This verifies part (b) of Condition (ii) Proposition  \ref{prop:AltFKG2} with
$p_0 = (1-p)^{2r_0} p > 0$.

For part (a) of Condition (ii) in Proposition \ref{prop:AltFKG2}, let us denote,
$I_{x} := h(x,0)(1) - x$, $I_{y} := h(y,0)(1) - y$
and $I_{z} := h(z,0)(1) - z$. As in the earlier model, we obtain,
\begin{align*}
& \E[X_{1}(1)X_{1}(2) - X_{0}(1)X_{0}(2) | X_0 = (y-x,z-y)]  \\
% = & \E \bigl[ ( y-x + (I_y - I_x))  (z-y + (I_z - I_y))  - (y-x)(z-y) | X_0 = (y-x,z-y) \bigr] \\
% = & \E \bigl[ ( y-x + (I_y - I_x))  (z-y + (I_z - I_y))  - (y-x)(z-y) \bigr] \\
& = (z-y) \E (I_y - I_x) +  (y-x) \E (I_z - I_y) + \E [ (I_y - I_x)(I_z - I_y) ].
\end{align*}
Now, we use the fact that the marginal distributions of $ I_x,I_y $
and $ I_z$ are same as $ I$ (say) with $\E(I) = 0$. But as opposed to the earlier model,
these random variables are no longer independent.
Using Cauchy Schwartz inequality for any $x , y , z \in \Z$ with $x < y < z$ we have
\begin{align}
\label{eq:CSbound}
& \E[X_{1}(1)X_{1}(2) - X_{0}(1)X_{0}(2) | X_0 = (y-x,z-y)] \nonumber \\
& =  \E [ (I_y - I_x)(I_z - I_y) ] \leq 4\E[I^2] < \infty.
\end{align}

To verify part (a) of the Condition (ii) of Proposition \ref{prop:AltFKG2},  we still need to show that
for $((y-x),(z-y)) \in (\mathcal{S}_0 \cup \mathcal{S}_1)^c$ there exists $d_1 > 0$ such that
$\E[X_{1}(1)X_{1}(2) - X_{0}(1)X_{0}(2) | X_0 = (y-x,z-y)] \leq -d_1$.
For $x < y < z \in \Z$ with $ \min (y -x, z-y)  > r > 0$ we have
\begin{align}
\label{eq:OutS1}
& \E[X_{1}(1)X_{1}(2) - X_{0}(1)X_{0}(2) | X_0 = (y-x,z-y)] \nonumber\\
& = \E[(I_y - I_x)(I_z - I_y)]  = - \E[I_y^2]  + \E [I_yI_z]  - \E[I_xI_z]  + \E[I_x I_y].
\end{align}
The last three terms converge to $0$ as $ r \to \infty $. We show it for $ \E [I_yI_z]  $, 
the others being exactly the same.  Observe that $ \E[ ( I_y I_z )\mathbb{I} 
(\max\{ |I_y|, |I_z|\} < r/2)] = 0$ and hence we obtain
\begin{align*}
\lefteqn{ |\E[ I_y I_z)] | } \\
& =  | \E[ ( I_y I_z )\mathbb{I} (\max\{|I_y|, |I_z|\} < r/2)] 
+ \E[ ( I_y I_z )\mathbb{I} (\max\{|I_y|, |I_z|\} \geq r/2)] | \\
& \leq \sqrt{\E[  I_y^2 I_z^2 ] } \sqrt{ \P (\max\{|I_y|, |I_z|\} \geq r/2)  } \\
& \leq \sqrt[4]{\E[  I_y^4] } \sqrt[4]{\E[  I_z^4] }  \sqrt{2} \sqrt{ \P ( |I_y| > r/2) } \\
& \to 0 
\end{align*}
as $ r \to \infty $. Thus, we have $ - \E[I_y^2]  + \E [I_yI_z]  - \E[I_xI_z]  + \E[I_x I_y] \leq - \E (I^2)/2 $
for all $ r > r_0 $. Hence it follows that Howard's model satisfies Proposition
\ref{prop:AltFKG2} with ${\cal S}_1 = \{(l,m) \in \Z^2 :  1 \leq \min \{ l, m \} \leq r_0 \}$, $d_1 =
  \E (I^2) /2, d_2 =  4\E[I^2] $ and $p_0 = (1-p)^{2r_0} p$
 and Condition (ii) of Proposition \ref{prop:AltFKG1} holds for this discrete model.
 This shows that the condition (B2) is satisfied.

\subsection{Ferrari, Landim and Thorisson's model}\label{sec:continuous}

We are interested in a directed forest spanning on Poisson process on $\R^2$
 introduced by Ferrari et al. \cite{FLT04}. The vertices of the forest
are the atoms of a homogeneous Poisson point process $N$ having unit intensity in $\R^2$, embedded
with the Euclidean norm. With an abuse of notation, the Poisson point process $N$
will be considered in the sequel sometimes as a point measure with Dirac masses
on its atoms and sometimes as a union of points of $\R^2$.
The ancestor of $(x,t)\in \R^2$, denoted by $h(x,t)$, is the closest point $(x', t')\in N$ such that
$$t'>t\qquad \mbox{ and }\qquad |x'-x|<\frac{1}{2}.$$
This means that in order to find the ancestor of $(x,t)$, we are growing a tube of
width 1 in the upward direction whose base is centered around $x$. The first point
of $N$ encountered by the tube is the ancestor of $(x,t)$.

The graph, considered by Ferrari et al. \cite{FLT04},  with vertices as the atoms of
$N$ and edge set $\big\{ \langle(x,t),h(x,t)\rangle,\ (x,t)\in N\big\}$, is of out-degree
$1$ almost surely and hence a forest. We should mention here that the authors
of \cite{FLT04} considered this forest on $\R^d$ for general $d$ and showed that
for $d=2$, the random forest is indeed a tree, i.e., connected  almost surely.
We concentrate only on $d=2$. In another work,  Ferrari, Fontes and Wu
\cite{FFW05}  showed that for $d=2$ under diffusive scaling, as a collection
of paths, this tree converges to the Brownian web.

%Condition (iii) of Proposition \ref{prop:AltFKG1} is not necessarily satisfied in the continuum case. 
It is not apriori clear that Condition (iii) of Proposition \ref{prop:AltFKG1} holds
in this model. However,
we show that the main idea of Proposition \ref{prop:AltFKG1} can be
extended to give an alternate verification of (B2) for this model as well.

For any path $ \pi \in {\cal X} $ emanating from $ (x,t)$, we define a c\`adl\`ag path, denoted
by $ \tilde{\pi} $ henceforth, emanating from $(x,t)$ to its ancestor $h(x,t) = (x',t')$
which remains constant on $[t,t')$ with value $x$ and then jumps to $x'$.
It is easy to observe that $ \pi $ and $ \tilde{\pi} $ meet  at jump times of $ \tilde{\pi} $ and
hence the number of steps  between two time segments are
exactly the same for both types of paths. From now on, we will only {\em consider
the c\`adl\`ag paths}.

In the sequel, we will deal with triplets $(x,y,z)\in \R^3$ 
with $x< y < z$ and $\min \{ (y-x), (z-y) \} \geq 1/2$.
 Let $X_t = ( \tilde{\pi}^{(x,0)} (t) , \tilde{\pi}^{(y,0)} (t) , \tilde{\pi}^{(z,0)} (t) ) $
denote the triplet of c\`adl\`ag paths starting from the
points $(x,0), (y, 0) $ and $ (z,0) $ respectively.  Further, it is useful to observe that
our choice of starting points ensures that
either $ \min \{ (X_t (2) - X_t (1) ), (X_t (3) - X_t (2) ) \} \geq 1/2$ or $ (X_t (2) - X_t (1) )
(X_t (3) - X_t (2) )  = 0 $ for all $ t \geq 0 $ where $ X_t (i) $ denotes the $i$-th co-ordinate
of $ X_t $. In other words,  the paths are separated by distance $1/2$ or they have coalesced
for all time points.
Hence, the process $ \{ X_t : t \geq 0 \}$,  is a pure jump Markov process with state space
$ {\cal M} = \{ (u,v, w) : u \leq v \leq w, \text{ either } u = v \text{ or }
v = w \text{ or } \min \{ v-u, w-v \} \geq 1/2 \}$.
%
%
%
%
% Further, properties of Poisson point
% process allows us to describe these c\`adl\`ag paths by Markov process, which
% we would exploit. It is useful to observe that, either two (c\`adl\`ag) paths  are separated by at least distance
% $1/2$  or they have coalesced.

% Next recall that Condition (iii) of Proposition \ref{prop:AltFKG2} does not
% hold for this model. Nevertheless, (\eqref{cond:exp-collision-time}) together
% with independence of Poisson process over disjoint regions, allow us to prove
% (B2) for this model.
This part of the argument is motivated from Ferrarri et. al.  \cite{FFW05}.
We are interested in finding $\lim_{n \to \infty} \P (\eta_{{\cal X}}(0,n;0,\epsilon \sqrt{n})
\geq 3)$. Set $A_n$ as the (random) point set $\{\tilde{\pi}(0) : \pi \in {\cal X}, \sigma_{\pi} \leq 0 \} \cap
 [-1/2,\epsilon\sqrt{n}+1/2]$. It is easy to observe that for any $x, y \in A_n$ we must have
$|x-y| \geq 1/2$. Hence elements of the set $A_n$ can be enumerated as
$\{x_1, x_2, \dotsc, x_J\}$ where $x_1, \dotsc, x_J$ and $ J$ are random
variables such that $x_1 < x_2 < \cdots < x_J$ and $ x_{i+1} - x_i \geq 1/2 $
almost surely for $ 1 \leq i \leq J-1 $ (with convention that if $ J = 0$, the point set is empty).
If required, we consider an enlarged set $B_n$ containing  $A_n$ such that
$B_n = \{y_1, y_2, \dotsc, y_{J^\prime}\}$ where
\begin{equation*}
y_1 \geq -1/2, y_{J^{\prime}} \leq \epsilon \sqrt{n} + 1/2 \text{ and } 1/2 \leq y_{i+1} - y_{i} \leq 1
\text{ for } 1 \leq i \leq J^{\prime} -1.
\end{equation*}
Note that $J^\prime$ is random but that we necessarily have $ J^\prime \leq 2\epsilon\sqrt{n} +2  $.

%\textcolor{blue}{
For  $1 \leq j \leq J^\prime - 2$, let us define the event
\begin{equation*}
E(j) = \{\tilde{\pi}^{(y_j,0)}(n) < \tilde{\pi}^{(y_{j + 1},0)}(n) < \tilde{\pi}^{(y_{J^\prime},0)}(n)\}.
\end{equation*}
On the event $\{\eta_{{\cal X}}(0,n;0,\epsilon \sqrt{n}) \geq 3\}$,
there must exist $j = j(\omega) \in \{1, 2, \cdots, J^\prime -2\}$ such that $E(j)$ is satisfied.
%\begin{equation*}
%\tilde{\pi}^{(y_j,0)}(n) < \tilde{\pi}^{(y_{j + 1},0)}(n) < \tilde{\pi}^{(y_{J^\prime},0)}(n).
%\end{equation*}
%}
Let ${\cal F}_0$ denote the $\sigma$-field generated by the Poisson
point process on the negative half-plane $\{(x,t)\in \R^2 : t \leq 0\}$.
Clearly, $A_n$ is ${\cal F}_0$ measurable. Given the set $A_n$, the set $ B_n $ can be chosen
by a deterministic algorithm so that $ B_n $, and hence $J^\prime$,  are ${\cal F}_0$ measurable.
On the other hand, given $B_n$, the event $E(j)$ depends only on
Poisson process on the upper half-plane $\{(x,t) \in \R^2 : t > 0\}$
and hence  $E(j)$ depends on ${\cal F}_0$ only through the positions of $ B_n $.

We prove the following proposition.
\begin{proposition}
\label{prop:FLTProbBound}
For any set $ B_n $ having properties as above, we have,  for $ j \geq 1 $,
\begin{equation*}
\P (E (j) \mid B_n) \leq \frac{1}{n} \bigl[ C_1 + C_2 \epsilon \sqrt{n} \bigr]
\end{equation*}
for constants $ C_1, C_2 > 0 $ independent of $ j $, as well as $B_n$.
\end{proposition}

Assuming the Proposition \ref{prop:FLTProbBound}, we have
\begin{align*}
& \P(\eta_{{\cal X}}(0,n;0,\epsilon \sqrt{n}) \geq 3)
= \E\bigl(\P((\eta_{{\cal X}}(0,n;0,\epsilon \sqrt{n}) \geq 3)| {\cal F}_0 )\bigr) \\
& \leq \E\bigl(\P(\cup_{j=1}^{J^\prime -2} E(j)| {\cal F}_0 )\bigr)
\leq \E\bigl(\sum_{j=1}^{J^\prime -2} \P( E(j)| {\cal F}_0 )\bigr) \\
& =  \E\bigl(\sum_{j=1}^{J^\prime -2} \P( E(j)| B_n )\bigr)
 \leq  \E\bigl(\sum_{j=1}^{J^\prime -2}\frac{1}{n} [C_1 + C_2\epsilon \sqrt{n}]\bigr) \\
& = \frac{1}{n}  [C_1 + C_2\epsilon \sqrt{n}] \E ( J^\prime -2) \leq \frac{1}{n}  [C_1 + C_2\epsilon \sqrt{n}]
2\epsilon\sqrt{n} \\
& \to 2 C_2 \epsilon^2
\end{align*}
as $ n \to \infty$, which proves (B2).

% The last step follows from (\eqref{cond:exp-collision-time}) and the observation that for all $1 \leq j \leq J^\prime - 2$, we have
% $ y_{j + 1} - y_j \leq 1 \text{ and } y_{J^\prime} - y_{j + 1}) \leq   \epsilon \sqrt{n} - j/2$.
% Finally $J^\prime \leq 2\epsilon\sqrt{n}$ gives that
% $$
% \P(\eta_{{\cal X}_1}(0,n;0,\epsilon \sqrt{n}) \geq 3) \leq  2C_1 \epsilon\sqrt{n} + C_3(\epsilon \sqrt{n})^2.
% $$
% Thus we have $\frac{1}{\epsilon}\limsup_{n \to \infty}\P(\eta_{{\cal X}_n}
%(0,n;0,\epsilon \sqrt{n}) \geq 3)\to 0$ as $\epsilon \downarrow 0$ and hence (B2) holds

\subsubsection{A Lyapunov condition for c\`adl\`ag branches}
% It is enough to show that the conditions of Proposition \ref{prop:AltFKG1} hold
% in this model.
% %Conditions (i), (ii), (iii) and (iv) (b) are easy to check if we choose
% %${\cal S}_1=\{(x,y)\in \R^2_+ : 0 < x\wedge y \leq x\vee y < 1\}$ and
% %${\cal S}_0=\{(x,y)\in \R^2_+ : xy = 0\}$.
% %%${\cal S}_1=\{(x,y,z) : 0 < (y-x)\wedge (z-y) \leq (y-x)\vee (z-y) < 1\}$ and
% %%${\cal S}_0=\{(x,y,z) : x \leq y \leq z, (y-x)(z-y) = 0\}$.
% Conditions (i) and (ii) are easy to check, hence we focus on the proof of Condition (iii).
%
% Here, instead of considering the forest
% we consider it as a collection of discontinuous paths. The discontinuous branch
% going through a point $(x,t)$ to its ancestor $h(x,t) = (x',t')$ remains constant
% on $[t,t')$ with value $x$ and then jump to $x'$.
% This allows us to describe the
% branches by Markov processes.
% We note here that the first collision time among
% three such discontinuous paths is same as that of the linearly interpolated
% paths.
%
%
%

It now remains to prove Proposition \ref{prop:FLTProbBound}. 
We are inspired by a result due to Meyn and Tweedie \cite[Theorem 4.3]{MT93-continuous},
which is based on establishing a Lyapunov control similar to Theorem
\ref{thm:AltFKGMarkov} (i) on the infinitesimal generator of $(X_t)_{t\geq 0}$
(see e.g. \cite{EK09} or \cite{MT93-continuous} for definitions of these generators).

Let us define $ {\cal M}_0 =  \{ (u,v, w) : u \leq v \leq w, \text{ either } u = v \text{ or }
v = w \} $ and $ \tau = \tau({\cal M}_0)=\inf\{t\geq 0,\ X_t \in {\cal M}_0\}$. We note that
$ \tau < \infty $ almost surely. 

\begin{theorem}\label{th:Lyapunov}
Assume that there exist $d > 0$ and  a non-negative function
$V\ :\ {\cal M} \rightarrow [0,\infty)$ such that the generator $G$ of the
Markov process $(X_t)_{t\geq 0}$ satisfies
\begin{equation}
\label{eq:Lyapunov-continuous}
GV(u,v, w)\leq - d  \ind_{ ({\cal M}_0)^c}(u,v, w) ,
\end{equation}
Then, for any $(u,v, w) \in  ({\cal M}_0)^c $,
$$\E_{(u,v, w)}(\tau)\leq \frac{1}{d} V(u,v, w) .$$
\end{theorem}

\begin{proof}[Theorem \ref{th:Lyapunov}] The result is a particular case of \cite{MT93-continuous} 
(Th. 4.3). Using Dynkin's formula:
\begin{align*}
0 & \leq \E_{(u,v, w)}\big(V(X_{\tau \wedge t})\big) \\
 & \leq  V(u,v, w) +\E_{(u,v, w)}\Big(\int_0^{\tau\wedge t}
\big(-d  \ind_{({\cal M}_0)^c}(X_s)\big) ds\Big)\\
& \leq  V\big(u,v, w\big) - d \E_{(u,v, w)}\big(\tau\wedge t\big)
\end{align*}
from which we obtain the announced result by letting $t \to \infty$.
\end{proof}

We claim that for the function 
\begin{equation}
\label{eqn:defV}
V\big(u,v, w\big) := (v-u)(w-v)
\end{equation}
the equation (\ref{eq:Lyapunov-continuous}) is satisfied with $ d = 1/12 $. Assuming the claim, 
we first prove Proposition \ref{prop:FLTProbBound}. 

\begin{proof}[Proposition \ref{prop:FLTProbBound}] With $ \tau $ being as defined before 
Theorem \ref{th:Lyapunov}, for $j\in \{1, \dotsc,  J'-2\}$, 
\begin{align*}
\P (E (j) \mid B_n) & = \P_{(y_j,y_{j+1},y_{J^\prime})}(\tau>n)
\leq \frac{1}{n}\E_{(y_j,y_{j+1},y_{J^\prime})}(\tau) \\
& \leq \frac{12}{n}   ( y_{j+1} - y_j ) (y_{J^\prime} -  y_{j+1})   \leq \frac{12}{n}    [\epsilon \sqrt{n} + 1]  
\end{align*}
using the fact that $ y_{j+1} - y_j \leq 1 $ and $ y_{J^\prime} -  y_{j+1} \leq \epsilon \sqrt{n} + 1 $. 
\end{proof}
% If we manage to upper bound the expectation $\E_{(y_j,y_{j+1},y_{J^\prime})}
%(\tau)$ by $C_1+C_2 \epsilon \sqrt{n}$, then \eqref{FFT-etape3} yields
% the inequality announced in Proposition \ref{prop:FLTExpBound}.}
% 
% 
% If we show that \eqref{eq:Lyapunov-continuous} is satisfied in our case,
% the theorem provides for $(u,v,w)=(y_j,y_{j+1},y_{J^\prime})$ that:
% \[\E_{(y_j,y_{j+1},y_{J^\prime})}(\tau)_leq \frac{(y_{J^\prime}-y_{j+1})
%(y_{j+1}-y_j)}{d}\leq \frac{\epsilon \sqrt{n}+1}{2d}\]
% which finishes the proof.

Now we verify the claim above.  For this, we start with computing the generator of $(X_t)_{t\geq 0}$.
The basic idea is the following: the generator computed on  $ V $ can be
interpreted as an expectation, which gives us a way to estimate it easily.
Let us fix $ (u,v, w) \in ({\cal M}_0)^c $.
We need some notation. Let $ {\cal J} = [-1/2, 1/2] $ and $ A =
(u+{\cal J}) \cup (v+{\cal J}) \cup (w+{\cal J})$. Note that these sets in $ A $
may overlap and the total length (Lebesgue measure) $|A|$ of $A$ can be strictly less than 3.
Consider a uniform random variable $U$ on the set $ A $. Now, for $s\in \{u,v, w\}$, define $I_s$ as follows:
\begin{equation*}
 I_s := \begin{cases} U - s & \text{ if } U \in s+{\mathcal J} \\ 0 & \text{ otherwise. } \end{cases}
\end{equation*}
Note that $I_s$ denotes the increment of the jump of the path starting from $(s,0)$ at the first jump time
of the process $ (X_t)_{t \geq 0 }$. Then we have
\begin{equation*}
G V (u,v, w) = | A |  \E \bigl[ V (u+I_u, v+I_v, w+I_w) - V (u,v, w) \bigr]. 
\end{equation*}

Using the form of $ V $, we have
\begin{equation}
\label{eqn:FLTG}
G V (u,v, w) = | A |  \E \bigl[  ( v-u ) (I_w - I_v) + (w-v) (I_v - I_u) + (I_w - I_v)  (I_v - I_u)  \bigr].
\end{equation}

First we make a few observations on $ I_s $ for $ s \in \{ u,v, w \} $. Clearly, $ I_s $ is symmetric
and hence $ \E ( I_s ) = 0 $.  Further, $ I_u I_w = 0 $ almost surely since $ (u+ {\cal J} ) \cap (w +{\cal J}) =
\emptyset $ as $ (u,v, w) \in ({\cal M}_0)^c $. Also, we observe that $ I_u I_v  $ is non-zero
only when $ U \in  (u+ {\cal J} ) \cap (v +{\cal J}) $. In such a case, $ I_u > 0 $ and $ I_v < 0 $. Thus,
$ I_u I_v  \leq 0$ almost surely. Similar argument  shows that $ I_v I_w  \leq 0$ almost surely.
Now, 
\begin{equation*}
\E ( I_v^2 ) = \frac{1}{|A|} \int_{-1/2}^{1/2}  t^2 dt = \frac{1}{ 12 | A| }.
\end{equation*}
Thus, putting  back in (\ref{eqn:FLTG}), we have
\begin{align*}
G V (u,v, w) & = | A |  \E \bigl[   (I_w - I_v)  (I_v - I_u)  \bigr] \leq -  |A| \E ( I_v^2 ) = - \frac{1}{ 12 }.
\end{align*}
This proves the claim that  $ V $ satisfies the condition (\ref{eq:Lyapunov-continuous}) with $ d = 1/12 $.

\section{Acknowledgement}
The authors thank Rongfeng Sun for careful reading and valuable comments. 
This work was supported by the GdR GeoSto 3477 and Labex CEMPI 
(ANR-11-LABX-0007-01). D.C. is funded by ANR PPPP (ANR-16-CE40-0016). 
A.S. wishes to thank Science \& Engineering Research Board (SERB), 
Department of Science \& Technology for the grant (MTR/2017/000293).
This work was partially done while K.S., A.S. and V.C.T. were visiting the 
Institute for Mathematical Sciences, National University of 
Singapore in 2017. The visit was supported by the Institute.

\end{document}